\documentclass[aos,preprint]{imsart}

\RequirePackage{amsthm,amsmath,amsfonts,amssymb,bbm}
\RequirePackage[authoryear]{natbib}
\RequirePackage{xcolor}
\RequirePackage[colorlinks,citecolor=blue,urlcolor=blue]{hyperref}
\RequirePackage{graphicx}
\RequirePackage{subfigure}
\RequirePackage{float}
\RequirePackage{multirow}
\RequirePackage{diagbox}
\RequirePackage[figuresright]{rotating}
\RequirePackage{pdflscape}
\RequirePackage{booktabs}
\RequirePackage{bm}
\RequirePackage{algorithm}
\RequirePackage{algorithmic}
\RequirePackage{changepage}
\RequirePackage{tikz}
\RequirePackage{enumerate}
\RequirePackage{lscape}
\RequirePackage{threeparttable}

\startlocaldefs

\theoremstyle{plain}
\newtheorem{theorem}{Theorem}[section]

\newtheorem{corollary}[theorem]{Corollary}

\theoremstyle{remark}
\newtheorem{assumption}{Assumption}[section]

\newtheorem{remark}{Remark}[section]

\newcommand{\ind}{\mathbbm{1}}

\newcounter{cnstcnt}

\newcounter{bnstcnt}

\makeatletter
\newcommand*{\centerfloat}{\parindent \z@
  \leftskip \z@ \@plus 1fil \@minus \marginparwidth
  \rightskip \leftskip
  \parfillskip \z@skip
}
\makeatother

\endlocaldefs

\begin{document}

\begin{frontmatter}
\pdfsubject{}

\title{Functional Principal Component Analysis: Minimax Optimality from Sparse to Dense Designs and Beyond} \runtitle{Functional PCA from Sparse to Dense Designs}

\begin{aug}
\author[A]{\fnms{Shunxing}~\snm{Yan}\ead[label=e1]{sxyan@tsinghua.edu.cn}}
\and
\author[B]{\fnms{Fang}~\snm{Yao}\ead[label=e2]{fyao@math.pku.edu.cn}\thanks{Author to whom correspondence may be addressed.}}

\address[A]{Department of Statistics and Data Science, Tsinghua University \printead[presep={ \ }]{e1}} \address[B]{School of Mathematical Sciences, Center for Statistical Science, Peking University  \printead[presep={ \ }]{e2}}
\end{aug}

\begin{abstract}
Functional data analysis is an important statistical field that treats data as random functions. Owing to its infinite-dimensionality, functional principal component analysis (FPCA) plays a central role as a dimension-reduction and representation tool. More importantly, the random functions are often not fully observed in reality, but instead measured with noise at discrete locations. While simpler problems, such as mean estimation, have been extensively studied for discretely observed data, a complete minimax theory for FPCA with diverging indices has remained one of the fundamental open problems for over two decades.  To fill this gap, we consider a sieve-based FPCA framework that pools discrete measurements across subjects and accommodates sampling frequencies ranging from sparse to dense under random designs.  We begin by establishing new minimax rates for covariance estimation over general bivariate function classes beyond tensor-product low-complexity structures, thereby identifying the intrinsically two-dimensional optimal rate. Our central theoretical contribution is to establish minimax optimal rates for both individually labelled eigenfunctions with diverging indices and leading principal subspaces of increasing dimensions, under explicit growth and regularity conditions. To establish these results, we develop a general framework for sharp minimax analysis of FPCA under discrete observation, separating spectral approximation from stochastic perturbations and exploiting their distinct structures, thereby offering a new approach to a broad class of nonparametric spectral estimation problems. For both spectral targets, the optimal rates reveal two distinct statistical mechanisms, spectral identification and nonparametric reconstruction, whose interplay gives rise to either four or two sampling regimes. This yields, for the first time, a precise characterization of the impact of discrete sampling across the full range from sparse to dense designs. Finally, as a downstream application, we analyze FPCA-based functional linear regression and show that the FPCA-based estimator attains the minimax-optimal rates, thereby resolving the long-standing question of whether the FPCA-based linear regression procedure can attain the optimality under discrete observations.
\end{abstract}



\end{frontmatter}

\section{Introduction}\label{sec:introduction}

Functional data analysis (FDA) has emerged as a standard framework for modeling data represented as functions over a domain, such as time or space;  see the monographs \citep{ramsay2005,ferraty2006nonparametric,hsing2015theoretical,kokoszka2017introduction} and the surveys \citep{wang2016functional,aneiros2019recent} for comprehensive treatments. Importantly, owing to the infinite-dimensionality of functional data, functional principal component analysis (FPCA) has become a fundamental dimension-reduction tool and has been applied to a wide range of downstream tasks  \citep{yao2005a,yao2005b,hall2007methodology,delaigle2012componentwise,dou2012estimation,chen2015localized, zhou2023functional}.  Despite substantial progress, a complete minimax theory for estimating eigenfunctions with diverging indices and leading principal subspaces of increasing dimension remains unavailable when, as is typically the case in practice, the trajectories are observed only at discrete locations. In this work, we address this central problem by developing minimax optimal procedures and theory for FPCA across the full range from sparse to dense observational designs.

Unlike the idealized setting of fully observed trajectories, functional data are typically observed only at discrete locations and are further contaminated by measurement error. Such discrete observations contain only limited information about each trajectory and introduce additional estimation errors. The resulting loss of information is particularly important for higher-order principal components, whose eigenvalues and eigengaps are typically shrinking towards zero, making the corresponding eigendirections more difficult to estimate. Nevertheless, estimating eigenfunctions with diverging indices and leading principal subspaces of increasing dimension is essential for many FPCA-based downstream analyses, where the number of retained components need to increase to control truncation bias. A theory restricted to fixed spectral indices therefore cannot characterize such procedures; the target spectral order must itself be allowed to increase with the sample size. Besides the stochasticity of underlying functions, the difficulty of estimating such growing spectral objects depends critically on the sampling frequency within each trajectory. Suppose that \(n\) independent trajectories are observed, with the \(i\)th trajectory measured at \(m_i\) locations. Assume that the deterministic counts satisfy \(m_i\geq2\) and \(m_i \asymp m\) uniformly over \(1 \leq i \leq n\), and refer to \(n\) and \(m\) as the sample size and sampling frequency, respectively.  Here and throughout, we write \( a_n \asymp b_n \) when both \( a_n \lesssim b_n \) and \( b_n \lesssim a_n \) hold, with  \( a_n \lesssim b_n \) (or \( a_n \gtrsim b_n \)) meaning that \( a_n \leq C b_n \) (or \( Ca_n \geq b_n \)) for some constant \( C > 0 \).  The sampling frequency \(m\) may remain bounded or increase with \(n\), thereby covering a continuum from sparse to dense observation designs. Hence it is important to study how the sample size \(n\), the sampling frequency \(m\), and the target spectral order \(j\) jointly determine the estimation error of FPCA in a unified framework, especially at what sampling frequency the fully observed statistical benchmark can be recovered.

\subsection{Related work and motivations}
The effect of discrete sampling on functional data estimation has been extensively studied, particularly for mean and covariance functions. Two main strategies have emerged. Under sufficiently dense designs, individual trajectories can be pre-smoothed and subsequently analyzed as if they were fully observed. This approach, however, may introduce substantial reconstruction error and is generally unsuitable for sparsely sampled data. By contrast, pooling observations across subjects avoids subject-specific reconstruction and provides a unified approach across sampling frequencies ranging from sparse to dense \citep{yao2005a,yao2005b,cai2011optimal,zhang2016sparse}. Existing studies have established convergence rates for mean and covariance estimation and revealed phase-transition phenomena as the sampling frequency increases \citep{cai2010nonparametric,cai2011optimal,li2010uniform, zhang2016sparse,xiao2018fast,xiao2020asymptotic,wang2022low,yan2025drrm}. Below a problem-dependent transition boundary, discrete sampling affects the convergence rate, which depends explicitly on the sampling frequency. Once the boundary is crossed, the fully observed rate is attained, and further increasing the number of observations per trajectory does not improve the convergence order. Nevertheless, existing minimax theory for covariance estimation has largely been developed under tensor spaces \citep{cai2010nonparametric,wang2022low}. Under such low-complexity structures, the bivariate covariance function behaves statistically like a one-dimensional smoothing problem, rather than a general two-dimensional function estimation problem. By contrast, \citet{li2010uniform}, \citet{zhang2016sparse} and \citet{xiao2020asymptotic} derive two-dimensional convergence rates, but without matching minimax lower bounds. Consequently, the minimax optimal rate of covariance estimation over general bivariate classes, without low-complexity tensor assumptions, remains unknown. A sharp covariance theory provides the necessary foundation for studying spectral recovery without imposing restrictive low-complexity structure. 

FPCA poses greater challenges than global covariance estimation because it requires recovering spectral structure whose stability deteriorates as the eigenvalues decay and the spectral gaps shrink. For fully observed trajectories, the optimal rate for estimating the \(j\)th eigenfunction is of order \(j^2/n\) \citep{hall2007methodology,wahl2022lower}. Although the adjacent eigengap is of order \(j^{-a-1}\), the variability of corresponding covariance perturbation also decreases as the corresponding eigenvalues decay, partially compensating for the shrinking gap. This yields the rate \(j^2/n\), rather than \(j^{2a+2}/n\) suggested by combining a global covariance bound with a Davis--Kahan inequality. Related refined perturbation analyses based on relative spectral separation and covariance-adapted perturbation bounds, have been developed \citep{JIRAK2020,JIRAK2023108808}. Under noisy discrete sampling, however, the convenient representation in terms of principal component scores available in the fully observed setting is no longer applicable. {Recovering the target eigendirection from discrete noisy measurements introduces an additional nonparametric reconstruction error, making the problem substantially more challenging. Previous studies have, for example, employed finite-dimensional basis representations and reduced-rank models \citep{james2000principal,peng2009geometric}, with convergence rates established for likelihood-based and penalized spline FPCA estimators \citep{paul2009consistency,he2024penalized}. Local smoothing methods have also been widely used \citep{yao2005a,hall2006properties,li2010uniform}. Nevertheless, existing theory for discretely observed FPCA has largely focused on a finite number of principal components \citep{cai2010nonparametric,10.1214/12-AOS1033,belhakem2025minimax}.} Recently, \citet{zhou2024theoryfunctionalprincipalcomponent} derived upper bounds for eigenfunctions with diverging indices from discrete observations.  Nevertheless, the available upper bounds are not sharp, especially for nonparametric reconstruction error of eigenfunctions, and lower bounds remain unavailable, leaving both the optimal rate and the precise phase-transition boundary unresolved. Existing results also focus primarily on individual eigenfunctions. For such targets, the rate-optimal tuning parameter generally depends on the eigenfunction index, so estimating multiple eigenfunctions at their individual optimal rates would require repeatedly constructing the covariance estimator with different tuning parameters. In practice, however, FPCA is typically based on a single covariance estimator, and subsequent analyses rely on the subspace spanned by its leading eigenfunctions. This motivates a direct study of leading principal subspace in increasing dimension, whose perturbation geometry differs fundamentally from that of individually labelled eigenfunctions.

There are many downstream tasks of FPCA, including trajectory recovery, regression, classification, and clustering, among which functional linear regression (FLR) is one of the most canonical models. As an ill-posed inverse problem, it is commonly regularized by truncating the principal components and estimating the slope function within the finite-dimensional leading principal subspace. Consistency requires the truncation level to diverge with the sample size, making a sharp FPCA theory for diverging spectral indices particularly important \citep{yao2005b,hall2007methodology,james2009functional,dou2012estimation}.  Partly because a sharp theory for FPCA under discrete sampling is unavailable, classical minimax theory for functional linear regression has largely been restricted to the idealized setting of fully observed trajectories \citep{hall2007methodology,cai2006prediction}. Optimal estimation under discrete observation schemes has remained unresolved for the past two decades \citep{yao2005b}. Recently, \citet{yanflr} addressed this problem by introducing a pooling-ridge estimator within the RKHS framework \citep{yuan2010reproducing,cai2012minimax} and establishing minimax bounds across varying sampling regimes.  When the eigenfunctions are assumed known, \citet{cai2026costdiscretizationfunctionallinear} also obtain the minimax results for FLR and propose data adaptive procedures for unknown eigenvalues. However, whether the FPCA-based FLR with unknown eigenfunctions estimated from the discretely observed data, can attain the same optimal rates of \citet{yanflr} remains unknown. Resolving this question requires sharp control of a growing collection of estimated eigendirections. Taken together, these gaps call for a minimax theory of FPCA that allows the spectral order to diverge and remains valid across the full range from sparse to dense observation designs.

\subsection{Main contributions}\label{sec:motcon}
Our main contributions are fourfold. At the beginning, we develop a sieve-based covariance estimator by pooling off-diagonal measurements across subjects, which directly accommodates sampling frequencies ranging from sparse to dense designs. The spectral decomposition of the estimated covariance operator yields a sieve-based FPCA estimator, whose rank-\(J\) truncation admits a reduced-rank representation \citep{james2000principal,paul2009consistency}. As a foundational result, we establish matching minimax upper and lower bounds for covariance estimation over general bivariate approximation classes without imposing tensor-product low-complexity structures. The resulting two-dimensional optimal risk is of order $ n^{-1}+(nm^2)^{-r/(r+1)}. $ To the best of our knowledge, this is the first result establishing the minimax optimality of the intrinsically two-dimensional rate for covariance estimation. It reveals a phase transition at \(m\asymp n^{1/(2r)}\): below this boundary, the risk is governed by the intrinsically two-dimensional nonparametric term \((nm^2)^{-r/(r+1)}\), whereas above it, the fully observed rate \(n^{-1}\) is attained and further increasing the sampling frequency yields no rate improvement.

\begin{figure}[t]
  \centering
  \begin{minipage}[t]{0.49\textwidth}
    \vspace{0pt}
    \centering
    \subfigure[Individually labelled eigendirection: \(p_1,p_2,p_3\) correspond, respectively, to the three transition scales in \eqref{eq:intro-pointwise-phase-boundaries}.
      \label{fig:pointwise-phase-diagram}]{\includegraphics[width=\linewidth]{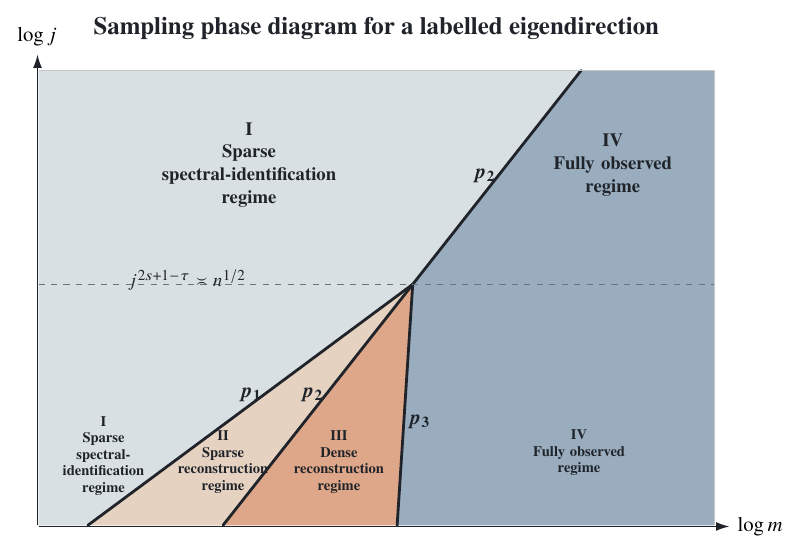}}
  \end{minipage}\hfill
  \begin{minipage}[t]{0.49\textwidth}
    \vspace{0pt}
    \centering
    \subfigure[Leading principal subspace: \(p_1,p_2,p_3\) correspond, respectively, to the three transition scales in \eqref{eq:intro-subspace-phase-boundaries}, up to logarithmic factors.
      \label{fig:subspace-phase-diagram}]{\includegraphics[width=\linewidth]{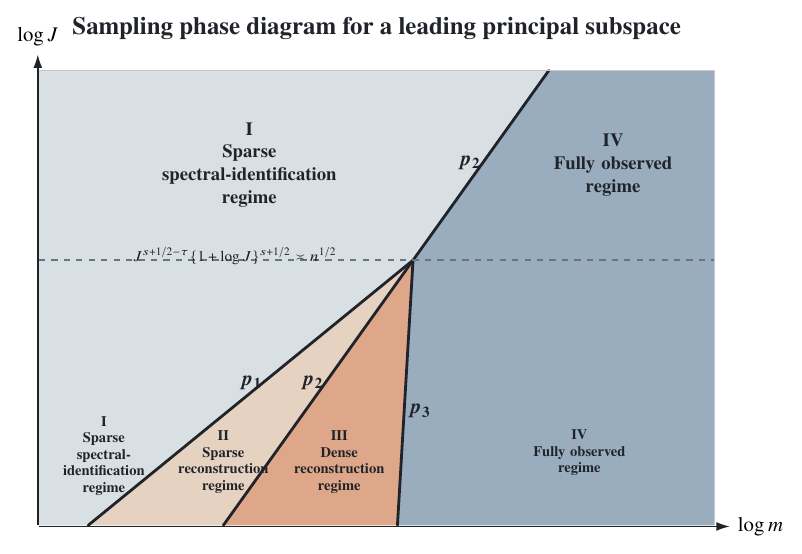}}
  \end{minipage}
  \caption{ Sparse-to-dense sampling phase diagrams for an individual  eigenfunction (left) and a leading principal subspace (right), plotted against the sampling frequency \(m\) and the target spectral order \(j\) or \(J\), respectively. The colored regions indicate the rate-determining regimes: light blue (I), sparse spectral identification; sand (II), sparse reconstruction; terracotta (III), dense reconstruction; and slate blue (IV), the fully observed benchmark. In each panel, the horizontal dashed line marks the critical spectral order separating the four-regime configuration below the line from the two-regime configuration above it. The displayed regimes occur under the parameter range described in Section~\ref{subsec:sampling-regimes}. }
  \label{fig:fpca-phase-diagrams}
\end{figure}

{ The central theoretical contribution is to establish matching minimax upper and lower bounds for an eigenfunction whose index \(j\) is allowed to diverge with the sample size.  To this end, we develop a general framework for sharp minimax analysis of FPCA under discrete observations, offering a new approach to a broad class of nonparametric spectral estimation problems. For the upper bound, we combine nonparametric approximation theory with refined spectral analysis, separating population spectral approximation from stochastic perturbations. The former is controlled at the intrinsic eigenfunction approximation rate without an additional eigengap loss, while the latter is analyzed by exploiting eigenvalue-dependent variances and individual spectral gaps in the eigenbasis of the projected covariance operator. For the lower bound, we provide two constructions, including a two-point rotation across a shrinking eigengap for spectral identification, and an Assouad hypercube of orthogonal rotations into a growing high-frequency block for nonparametric reconstruction; which capture the two distinct mechanisms underlying the minimax rate. } On the admissible domains specified by the corresponding bounds in Section~\ref{sec:theory}, the optimal rate in terms of squared \(\mathcal{L}^2\)-loss is
\[
\underbrace{
\frac{j^2}{n}
\left(1+\frac{j^a}{m}\right)^2
}_{\text{local spectral identification}}
+ \ 
\underbrace{
\frac{j^{(2as+2\tau)/(2s+1)}}
{(nm)^{2s/(2s+1)}}
\left(1+\frac{j^a}{m}\right)^{2s/(2s+1)}
}_{\text{nonparametric reconstruction error}}.
\]
Here, \(a\) characterizes the polynomial eigenvalue decay, while \(s\) and \(\tau\) govern the sieve approximation complexity, with \(j^\tau K^{-s}\) capturing the index-dependent approximation error of the \(j\)th eigendirection at sieve resolution \(K\). The two terms reflect distinct statistical challenges. The first arises from local spectral direction identification under shrinking eigengaps, whereas the second is a nonparametric reconstruction error associated with recovering the target eigendirection from discrete noisy observations. Their competition gives rise to a refined sparse-to-dense phase diagram over the parameter range described in Section~\ref{subsec:sampling-regimes}; see Figure~\ref{fig:pointwise-phase-diagram}. Inspecting the lower half of the plot when $ j^{2s+1-\tau}=o(n^{1/2}), $  the rate admits a four-regime characterization, separated by the transition scales
\begin{equation}\label{eq:intro-pointwise-phase-boundaries}
m\asymp n^{-1/2}j^{a+2s+1-\tau},
\qquad
m\asymp j^a,
\qquad
m\asymp n^{1/(2s)}j^{a-2+(\tau-1)/s},
\end{equation}
corresponding to $p_1, p_2$ and $p_3$ labelled in Figure~\ref{fig:pointwise-phase-diagram}, respectively. Here, the regimes are understood over their intersections with the admissible range. {The first and third boundaries arise from the competition between local spectral identification and nonparametric reconstruction. The middle boundary \(m\asymp j^a\) has a direct statistical interpretation, since \(j^a/m\asymp (m\mu_j)^{-1}\) is the inverse within-trajectory signal-to-noise ratio for the \(j\)th principal component, where $\mu_j$ is the $j$th nonzero eigenvalue of the covariance.  As \(m\) increases across these boundaries, the dominant source of error first shifts from sparse local spectral identification to sparse nonparametric reconstruction. At the second boundary, the \(j\)th principal component score becomes sufficiently well resolved within each trajectory, while nonparametric reconstruction remains rate-determining. Beyond the third boundary, the reconstruction error under dense sampling becomes negligible, and local spectral identification again determines the rate, which reaches the fully observed rate \(j^2/n\). Intuitively, a relative small eigengap makes the eigendirection harder to distinguish from its neighbors. When \(j^{2s+1-\tau}\gtrsim n^{1/2}\), the spectral-identification term remains rate-determining throughout the sampling range. Consequently, nonparametric reconstruction no longer gives rise to separate regimes, leaving only the transition around \(m\asymp j^a\) between sparse and dense spectral identification.}

Moreover, we establish a separate minimax theory for estimating the leading principal subspace of increasing dimension \(J\). Unlike individually labelled eigenfunctions, a principal subspace is invariant under rotations within its range. Consequently, simply aggregating individual eigenfunction errors does not yield a sharp characterization of the subspace estimation error. Under the normalized squared Hilbert--Schmidt loss, we establish the optimal rate
\[
\underbrace{
\frac{J(1+\log J)}{n}
\left(1+\frac{J^a}{m}\right)^2
}_{\text{cross-boundary spectral identification}}
+ \ 
\underbrace{
\frac{
J^{(2as+2\tau)/(2s+1)}
}{
(nm)^{2s/(2s+1)}
}
\left(1+\frac{J^a}{m}\right)^{2s/(2s+1)}
}_{\text{nonparametric reconstruction error}}.
\]
The first term is fundamentally different from its elementwise counterpart. Since rotations within the leading subspace do not change its span, only perturbations crossing the spectral boundary between the first \(J\) eigendirections and their orthogonal complement contribute to the loss. A direct analysis of the spectral projector exploits this geometry and yields a spectral-identification error of order \(n^{-1}J(1+\log J)(1+J^a/m)^2\), improving upon the \(n^{-1}J^2(1+J^a/m)^2\) order obtained by aggregating elementwise eigenfunction errors. In contrast, the nonparametric reconstruction error remains of the same order as the corresponding term for the eigenfunction \(\phi_J\). Thus, eigenfunction and subspace estimation share the same nonparametric reconstruction error, whereas their spectral-identification errors differ because only perturbations crossing the boundary of the leading subspace affect the spectral projector. The interplay between these two mechanisms gives rise to an analogous phase diagram for leading principal subspace estimation; see Figure~\ref{fig:subspace-phase-diagram}. When $ J^{s+1/2-\tau}(1+\log J)^{s+1/2}=o(n^{1/2}), $ the rate admits a four-regime characterization, separated by
\begin{equation}\label{eq:intro-subspace-phase-boundaries}
m\asymp
n^{-1/2}J^{a+s+1/2-\tau},
\qquad
m\asymp J^a,
\qquad
m\asymp
n^{1/(2s)}
J^{a-1+(\tau-1/2)/s},
\end{equation}
up to logarithmic factors, where the regimes are understood over their intersections with the admissible sampling range. By contrast, when \(J^{s+1/2-\tau}(1+\log J)^{s+1/2}/n^{1/2}\to\infty\), cross-boundary spectral identification dominates and only the two regimes separated around \(m\asymp J^a\) remain.  Importantly, because spectral identification is less difficult for a principal subspace than for the labelled eigenfunction \(\phi_J\), the four-regime structure persists over a substantially broader range of spectral dimensions. 

As a downstream application of the spectral analysis, we consider FPCA-based functional linear regression with discretely observed predictors. With a diverging truncation level, we show that the FPCA-based estimator attains the minimax-optimal rates established in \citet{yanflr} for both the squared \(\mathcal{L}^2\)-loss $ \|\widehat\beta-\beta\|_{\mathcal{L}^2(\mathcal T)}^2 =O_p\big( n^{-(2b-1)/(a+2b)} + (nm)^{-(2b-1)/(2a+2b)} \big) $, and the prediction loss $ \|C^{1/2}(\widehat\beta-\beta)\|_{\mathcal{L}^2(\mathcal T)}^2 =O_p\big( n^{-(a+2b-1)/(a+2b)} + (nm)^{-(a+2b-1)/(2a+2b)}\big), $ respectively. Here, \(b\) characterizes the regularity of the slope function through the decay of its coefficients in the covariance eigenbasis. We also locate the phase-transition boundary at $ m\asymp n^{a/(a+2b)} $.  This establishes the minimax optimality of the widely used FPCA-based regression procedure across sampling regimes and demonstrates how sharp FPCA theory for diverging spectral indices translates to a canonical functional inverse problem.

The paper is organized as follows. In Section~\ref{sec:fpca}, after describing the discrete observation model, we introduce the sieve-based covariance estimator and the FPCA counterpart, and establish the corresponding minimax-optimal bivariate rate for covariance estimation. Building on these, Section~\ref{sec:theory} develops matching minimax upper and lower bounds for estimating eigenfunctions with diverging indices and leading principal subspaces of increasing dimension, while delineating the associated sparse-to-dense phase transitions. Section~\ref{sec:functional-linear-regression} then turns to FPCA-based functional linear regression and establishes its minimax optimality under discrete sampling.  All proofs and auxiliary lemmas are deferred to the supplementary material.

\section{Discretely Observed Functional Data and Principal Component Analysis}\label{sec:fpca}
We begin with the standard framework for functional data. Let \(X\) be a square-integrable random function on the compact domain \(\mathcal T=[0,1]\), and let \(X_1,\ldots,X_n\) be independent and identically distributed copies of \(X\). Its mean function is $ \mu_X(t)=\mathbb E\{X(t)\}. $ Since mean estimation for discretely observed functional data is well studied \citep{cai2011optimal,zhang2016sparse,yan2025drrm}, we focus on the covariance and spectral structure and assume throughout that the mean function is zero on all \(\mathcal T\) for simplicity.  The covariance function of \(X\) is defined by $ C(t_1,t_2) = \mathbb E\left[  \{X(t_1)-\mu_X(t_1)\} \{X(t_2)-\mu_X(t_2)\} \right], $ which is symmetric and nonnegative definite. In the following, we assume that $\sup_{t\in\mathcal T} \mathbb{E} |X(t)|^4 \lesssim 1$.  Under the standard continuity conditions, Mercer's theorem gives the following spectral decomposition
\[
    C(t_1,t_2)
    =
    \sum_{k=1}^{\infty}
    \mu_k\phi_k(t_1)\phi_k(t_2),
\]
where \(\mu_1\geq\mu_2\geq\cdots>0\) are the nonzero eigenvalues of the covariance operator and \(\{\phi_k:k\geq1\}\) are the corresponding orthonormal eigenfunctions in \(\mathcal{L}^2(\mathcal T)\). Accordingly, each trajectory \(X_i\) admits the Karhunen--Lo\`eve expansion
\[
    X_i(t)
    =
    \sum_{k=1}^{\infty}\xi_{ik}\phi_k(t),
    \qquad
    \xi_{ik}
    =
    \int_{\mathcal T}X_i(t)\phi_k(t)\,\mathrm{d}t,
\]
where \(\xi_{ik}\) is the \(k\)th functional principal component score of $X_i$, with $ \mathbb E(\xi_{ik})=0 $ and $ \mathbb E(\xi_{ik}\xi_{i\ell})=\mu_k\ind\{k=\ell\}. $

In practice, the complete trajectories \(X_1,\ldots,X_n\) are not observed. Instead, the \(i\)th trajectory is measured with noise at \(m_i\) discrete locations: $ \bigl\{(T_{ij},X_{ij}):1\leq j\leq m_i\bigr\}, $ where \(T_{ij}\in\mathcal T\) is the \(j\)th observation time for subject \(i\). In the random-design setting, we assume that the observation times \(T_{ij}\)'s are i.i.d. sampled from a uniform distribution on \(\mathcal T\) and independent of other random variables. The observations follow the model
\begin{equation}
    X_{ij}
    =
    X_i(T_{ij})+\varepsilon_{ij},
    \qquad
    1\leq i\leq n,\quad
    1\leq j\leq m_i,
    \label{equ:lisanx}
\end{equation}
where \(\{\varepsilon_{ij}\}\) are i.i.d. sub-Gaussian measurement errors with mean zero and variance $\sigma_\varepsilon^2$, independent of all other random variables. For simplicity, we assume that \(m_i\geq2\) and \(m_i\asymp m\) uniformly over \(1\leq i\leq n\), and use \(m\) as the common notation of the subject-specific sampling frequencies. Observations from different subjects are independent, whereas observations from the same subject are dependent through the shared trajectory \(X_i\). 

\subsection{Sieve covariance and FPCA estimator}\label{sec:fpca_estimator}

To estimate the covariance operator from the discretely observed data, we introduce a fixed sequence of nested sieve spaces
\[
    \mathcal V_0=\{0\}
    \subset \mathcal V_1
    \subset \mathcal V_2
    \subset \cdots
    \subset \mathcal{L}^2(\mathcal T),
    \qquad
    \dim(\mathcal V_K)=K.
\]
These spaces may be generated by standard approximation systems, such as splines, Fourier bases, or wavelets, and are assumed to satisfy $ \overline{\bigcup_{K\geq1}\mathcal V_K} = \mathcal{L}^2(\mathcal T). $ Let \(\Pi_K\) denote the \(\mathcal{L}^2(\mathcal T)\)-orthogonal projector onto \(\mathcal V_K\).  For each \(t\in\mathcal T\), let \(g_{K,t}\in\mathcal V_K\) denote the projection representer satisfying
\[
    (\Pi_K f)(t)
    =
    \langle f,g_{K,t}\rangle_{\mathcal L^2},
    \qquad
    f\in\mathcal L^2(\mathcal T).
\]
The function \(g_{K,t}\) depends only on the sieve space \(\mathcal V_K\) and the point \(t\), and is invariant to the particular basis used to represent \(\mathcal V_K\). For implementation, we may represent the sieve space \(\mathcal V_K\) by a linearly independent basis \( \boldsymbol b_K(t) = \bigl( b_{K1}(t),\ldots,b_{KK}(t) \bigr)^\top \), and define its Gram matrix by $ G_K = \int_{\mathcal T} \boldsymbol b_K(t)\boldsymbol b_K(t)^\top \,\mathrm{d}t, $ then we have $ g_{K,t}(\cdot) = \boldsymbol b_K(\cdot)^\top G_K^{-1} \boldsymbol b_K(t). $ 

At a given sieve dimension \(K\), we define the off-diagonal covariance estimator by
\begin{equation}
    \widehat C
    =
    \frac{1}{n}
    \sum_{i=1}^n
    \frac{1}{m_i(m_i-1)}
    \sum_{j\neq\ell}
    X_{ij}X_{i\ell}
    g_{K,T_{ij}}\otimes g_{K,T_{i\ell}},
    \label{equ:covariance_estimator}
\end{equation}
where \(u\otimes v\) denotes the rank-one operator $ (u\otimes v)f = \langle f,v\rangle u. $ The use of off-diagonal products removes the bias induced by measurement errors. Indeed, for \(j\neq \ell\), the independence and mean-zero assumptions on the measurement errors imply $ \mathbb E\!\left( X_{ij}X_{i\ell} \mid T_{ij},T_{i\ell} \right) = C(T_{ij},T_{i\ell}). $ Together with the random design, this yields
\[
    \mathbb E(\widehat C)
    =
    \Pi_K C\Pi_K.
\]
Thus, \(\widehat C\) is an unbiased estimator of the covariance operator projected onto the sieve space.

Since \(\widehat C\) is a self-adjoint finite-rank operator on \(\mathcal V_K\), it admits the spectral decomposition
\[
    \widehat C
    =
    \sum_{\ell=1}^{K}
    \widehat\mu_{\ell}
    \widehat\phi_{\ell}
    \otimes
    \widehat\phi_{\ell},
\]
where the eigenvalues are arranged in decreasing order and \(\{\widehat\phi_{\ell}:1\leq\ell\leq K\}\) are the corresponding orthonormal eigenfunctions. Without loss of generality, we choose their signs so that $ \langle \widehat\phi_\ell,\phi_\ell\rangle_{\mathcal{L}^2(\mathcal T)} \geq 0 $ for all \(\ell \geq 1\). For \(J\leq K\), define the population and empirical projectors onto the leading \(J\)-dimensional principal subspaces by
\begin{equation}\label{equ:yinruPJ}
    Q_J
    =
    \sum_{j=1}^J
    \phi_j\otimes\phi_j,
    \qquad
    \widehat Q_{J}
    =
    \sum_{j=1}^J
    \widehat\phi_{j}\otimes\widehat\phi_{j},
\end{equation}
Then \(\widehat Q_{J}\) provides a natural estimator of the leading principal subspace.

The estimator above follows a standard sieve-based FPCA construction: it first estimates the covariance operator within the sieve space and then extracts its leading eigenfunctions. Its rank-\(J\) truncation also admits the familiar reduced-rank representation used in functional principal component models \citep{james2000principal,paul2009consistency}. Specifically, define the rank-\(J\) truncation $ \widehat C_J = \sum_{\ell=1}^{J} \widehat\mu_\ell \widehat\phi_\ell \otimes \widehat\phi_\ell, $ which retains the first $J$ empirical eigendirections of \(\widehat C\). In the basis representation \(\widehat\phi_\ell(t)=\boldsymbol b_K(t)^\top\widehat c_\ell\), let $ \widehat\Theta_J = (\widehat c_1,\ldots,\widehat c_J), $ and $ \widehat D_J = \operatorname{diag} (\widehat\mu_1,\ldots,\widehat\mu_J). $ Then
\[
    \widehat C_J(s,t)
    =
    \boldsymbol b_K(s)^\top
    \widehat\Theta_J
    \widehat D_J
    \widehat\Theta_J^\top
    \boldsymbol b_K(t),
    \qquad
    \widehat\Theta_J^\top G_K\widehat\Theta_J=I_J.
\]
Thus, \(\widehat C_J\) admits the usual reduced-rank sieve representation. When \(\mathcal V_K\) is generated by splines, this reduces to the familiar reduced-rank spline representation.

\subsection{Minimax Covariance Estimation over General Bivariate Approximation Classes}
\label{sec:covariance_minimax}

Existing sharp minimax theory for covariance estimation from discretely observed functional data has largely focused on tensor-product or related low-complexity function classes. In particular, \citet{cai2010nonparametric} model the covariance function in a tensor-product RKHS and obtain rates corresponding to essentially one-dimensional smoothing. Another line of work estimates the covariance function by directly smoothing off-diagonal raw products. Under more general bivariate smoothness conditions, \citet{li2010uniform,zhang2016sparse} derive convergence rates over sampling frequencies ranging from sparse to dense designs. Matching minimax lower bounds, however, are not available in these settings. Consequently, the minimax rate for covariance estimation over general bivariate approximation classes remains unknown.

We first impose a standard condition that controls point evaluation on the sieve spaces. 
\begin{assumption}[Sieve Point-evaluation Stability]
\label{assump:sieve_stability}
For the sieve spaces \(\mathcal V_K\),
\[
    \|v\|_{\mathcal L^\infty(\mathcal T)}
    \lesssim
    K^{1/2}\|v\|_{\mathcal L^2(\mathcal T)}
\]
uniformly over \(K\geq1\) and \(v\in\mathcal V_K\).
\end{assumption}

Equivalently, Assumption~\ref{assump:sieve_stability} requires $ \sup_{t\in\mathcal T} \|g_{K,t}\|_{\mathcal{L}^2(\mathcal T)} \lesssim K^{1/2}, $ which is the norm of the point-evaluation functional on \(\mathcal V_K\), and thus quantifies the stability of the sieve.  This condition controls the pointwise concentration of the sieve and prevents any single observation from having an unduly large effect after projection onto \(\mathcal V_K\). Obviously, this condition is satisfied by trigonometric sieves and other sieves generated by uniformly bounded orthonormal systems. More generally, it also holds for histogram and fixed-degree spline spaces on quasi-uniform partitions, as well as for compactly supported boundary-corrected wavelet spaces under their natural multiresolution ordering \citep{stone1994use,huang2003local,brenner2008mathematical,ginebook}.

The following condition is for the smoothness of the covariance function.
\begin{assumption}[Covariance Smoothness]\label{assump:covariance_smoothness}
For a fixed \(r>1/2\),
$$
    \|C-\Pi_KC\Pi_K\|_{\mathrm{HS}}^2
    \lesssim K^{-2r},
$$
where \(\Pi_K\) is the \(\mathcal{L}^2(\mathcal T)\)-orthogonal projector onto the sieve space \(\mathcal V_K\) defined in Section~\ref{sec:fpca_estimator}.
\end{assumption}

Assumption~\ref{assump:covariance_smoothness} is a sieve approximation condition and is commonly satisfied by a wide range of smooth function classes.  Under standard approximation conditions for spline or wavelet sieves, it corresponds to a Besov  \(B^r_{2,\infty}(\mathcal T^2)\)-type approximation condition. For periodic Fourier sieves, it requires polynomial decay of the high-frequency Fourier tail; in particular, periodic Sobolev smoothness of order \(r\) is sufficient.

The condition is weaker than the commonly used tensor-product RKHS condition \citep{cai2010nonparametric}. To illustrate this, let \(G\) be a Mercer kernel with eigenvalues \(\eta_k\asymp k^{-2r}\) and eigenfunctions \(\psi_k\), and take \(\mathcal V_K=\operatorname{span}\{\psi_1,\ldots,\psi_K\}\). Consider the covariance operator $ C = \sum_{k\geq1} k^{-r-1/2} \psi_k\otimes\psi_k. $ Then $ \|C-\Pi_KC\Pi_K\|_{\mathrm{HS}}^2 \lesssim K^{-2r}, $ so Assumption~\ref{assump:covariance_smoothness} holds. However, the corresponding covariance kernel does not belong to the tensor-product RKHS \(\mathcal H_G\otimes\mathcal H_G\), since $ \sum_{k\geq1} {(k^{-r-1/2})^2}/{\eta_k^2} = \infty $ for \(r>1/2\). Thus, the proposed condition includes more functions than the usual tensor RKHS class with smoothness $r$. As shown below, enlarging the covariance class changes the nonparametric minimax rate from the essentially one-dimensional rate induced by tensor-product RKHS regularity to an intrinsically two-dimensional rate.

\begin{theorem}[Lower Bound for Covariance Estimation]
\label{thm:covariance_lower}
Let \(\mathcal Q_{\mathrm{cov}}\) denote the collection of all data distributions satisfying Assumptions~\ref{assump:sieve_stability}--\ref{assump:covariance_smoothness} with a common fixed constant. Then there exists a constant \(c>0\), independent of \(n\) and \(m\), such that
\[
\varliminf_{n\to\infty}
    \inf_{\widetilde C}
    \sup_{\mathcal P\in\mathcal Q_{\mathrm{cov}}}
    \mathbb P_{\mathcal P}\left\{
    \bigl\|
        \widetilde C-C_{\mathcal P}
    \bigr\|_{\mathrm{HS}}^2
    \geq
    c\left[
        n^{-1}
        +(nm^2)^{-r/(r+1)}
    \right]
    \right\}
    >0,
\]
where the infimum is taken over all covariance estimators \(\widetilde C\) based on the observed data $ \bigl\{ (T_{ij},X_{ij}): 1\leq i\leq n,\; 1\leq j\leq m_i \bigr\}, $ and \(C_{\mathcal P}\) denotes the covariance operator under \(\mathcal P\).
\end{theorem}

The term \(n^{-1}\) persists even when the trajectories are fully observed and therefore represents the irreducible error from \(n\) independent subjects. The second term, $ (nm^2)^{-r/(r+1)} = (nm^2)^{-2r/(2r+2)}, $ arises from estimating a block of covariance coefficients indexed by two varying sieve coordinates. Because the covariance class does not impose product-type decay within this block, a full two-dimensional array of coefficients remains statistically relevant. The resulting nonparametric rate is therefore that of estimating a general bivariate smooth function. This differs fundamentally from the tensor-product setting of \citet{cai2010nonparametric}, where the imposed joint regularity effectively reduces the nonparametric estimation problem to essentially one-dimensional smoothing and yields the rate $ n^{-1}+(nm)^{-2r/(2r+1)}. $ In contrast, the present covariance class gives the intrinsically two-dimensional rate $ n^{-1}+(nm^2)^{-2r/(2r+2)}. $ The following upper bound shows that the proposed sieve covariance estimator attains this minimax bound and establishes optimality.

\begin{theorem}[Upper Bound for Covariance Estimation]
\label{thm:covariance_upper}
Under Assumptions~\ref{assump:sieve_stability}--\ref{assump:covariance_smoothness}, for any sieve dimension \(K\), the covariance estimator in \eqref{equ:covariance_estimator} satisfies
\[
    \|\widehat C-C\|_{\mathrm{HS}}^2
    = O_p\left(
    K^{-2r}
    +
    \frac{1}{n}
    +
    \frac{K}{nm}
    +
    \frac{K^2}{nm^2}\right).
\]
In particular, choosing a positive integer $ K \asymp (nm^2)^{1/(2r+2)} $ yields
\[
    \|\widehat C-C\|_{\mathrm{HS}}^2
    = O_p\left(
    n^{-1}+(nm^2)^{-r/(r+1)}\right).
\]
\end{theorem}

The four terms in Theorem~\ref{thm:covariance_upper} before trade-off have distinct origins. The term \(K^{-2r}\) is the sieve approximation error, while \(n^{-1}\) is the subject-level error from the \(n\) independent trajectories. The term \(K/(nm)\) arises from the interaction between the projected trajectory and a single discretization error and therefore scales with \(K\) sieve coordinates. In contrast, \(K^2/(nm^2)\) arises when discretization errors enter both arguments of the covariance operator and therefore scales with a full \(K\times K\) coefficient block. Note that \(K/(nm) \lesssim n^{-1} + K^2/(nm^2)\), so the \(K/(nm)\) term never determines the overall rate. Balancing \(K^{-2r}\) and \(K^2/(nm^2)\) gives \(K\asymp(nm^2)^{1/(2r+2)}\), yielding the minimax-optimal rate in Theorem~\ref{thm:covariance_upper}.

Finally, the minimax rate exhibits a phase transition at \(m\asymp n^{1/(2r)}\). In the sparse regime \(m\lesssim n^{1/(2r)}\), discrete sampling is rate-determining and the convergence rate is of order \((nm^2)^{-r/(r+1)}\). In the dense regime \(m\gtrsim n^{1/(2r)}\), the fully observed benchmark \(n^{-1}\) is attained, and further increases in the sampling frequency do not improve the rate.

\section{Minimax Theory for Eigenfunctions and Principal Subspaces}\label{sec:theory}

Although the covariance estimator studied in Section~\ref{sec:covariance_minimax} is minimax optimal under the Hilbert--Schmidt loss, the global error bound is too coarse for sharp spectral analysis at growing indices, where shrinking eigengaps amplify local perturbations. In this section, we develop minimax theory for individual eigenfunctions and leading principal subspaces whose spectral order is allowed to increase with the sample size. The following assumption is standard for functional data.

\begin{assumption}[Concentration and Spectrum]\label{ass:cov} Consider the following conditions. 
\begin{enumerate}[(i)]
\item  Suppose that \(\sup_{t \in \mathcal T}\left\|X(t)\right\|_{\psi_2} \lesssim 1\) and that, uniformly over \(f\in \mathcal{L}^2(\mathcal T)\), \[
    \big\|\left\langle X,f\right\rangle_{ \mathcal{L}^2(\mathcal T)}\big\|_{\psi_2}^2
    \lesssim
    \left\langle f,Cf\right\rangle_{ \mathcal{L}^2(\mathcal T)}.
\]
\item There exists a constant \(a>1\) such that, uniformly over \(j\geq1\),
\[
    \mu_j\asymp j^{-a},
    \qquad
    \mu_j-\mu_{j+1}\gtrsim j^{-a-1}.
\]

\item  There exist constants \(\tau,s>0\) and \( e \geq 1\), uniformly over $j,K$ with $1 \leq j \leq (K/2)^{1/e}$,
\[
    \|(I-\Pi_K)\phi_j\|_{\mathcal L^2(\mathcal T)}
    \lesssim j^\tau K^{-s}.
\]
\end{enumerate}
\end{assumption}

Assumption~\ref{ass:cov} collects the main concentration and spectral conditions on \(X\). Condition~(i) imposes sub-Gaussian concentration on both pointwise observations and linear functionals of \(X\). Condition~(ii) specifies polynomial eigenvalue decay and separation between adjacent eigenvalues. It determines the stability of spectral estimation and the degree of ill-posedness in associated inverse problems \citep{hall2007methodology,zhou2024theoryfunctionalprincipalcomponent}.

Condition~(iii) is on the smoothness regularity of the functions. Since the functional data are observed discretely, smoothness is necessary for eigenfunction recovery. For example, for Fourier-type sieves, a convenient sufficient condition is $ \|\phi_j\|_{H^s(\mathcal T)} \lesssim j^\tau $ under the usual boundary conditions associated with the chosen Fourier basis. Standard Fourier approximation then gives $ \|(I-\Pi_K)\phi_j\|_{\mathcal L^2(\mathcal T)} \lesssim K^{-s}\|\phi_j\|_{H^s(\mathcal T)} \lesssim j^\tau K^{-s}. $ Analogous bounds hold for spline and wavelet sieves by standard Jackson-type approximation inequalities, with the corresponding approximation classes naturally related to Besov \(B^s_{2,\infty}(\mathcal T)\)-type smoothness. In previous functional data works \citep{hall2006properties,zhou2024theoryfunctionalprincipalcomponent}, second-order smoothness with \(s=2\) is often considered, while the derivative-growth conditions imposed by \citet{zhou2024theoryfunctionalprincipalcomponent} suggest growth exponent \(2\) for common Fourier and wavelet systems, corresponding to \(\tau=2\) here.

Condition~(iii) plays an important role for the bias of the  eigenfunction estimation in the subsequent analysis. Consider the finite-dimensional projected covariance operator $ C_K=\Pi_K C\Pi_K . $ For \(1\leq j\leq K\), let \(\phi_{j,K}\) denote a unit eigenfunction associated with the \(j\)th ordered eigenvalue of \(C_K\). If this eigenvalue has multiplicity greater than one, we choose \(\phi_{j,K}\) within the corresponding eigenspace to maximize $ \bigl| \langle \phi_{j,K},\phi_j\rangle_{\mathcal L^2(\mathcal T)} \bigr|. $ To quantify how well the sieve preserves the \(j\)th eigendirection, define \(\Theta_{j,K}(C)\in[0,\pi/2]\) as the angle between \(\phi_j\) and \(\phi_{j,K}\), then
\[
    \cos\Theta_{j,K}(C)
    =
    \bigl|
        \langle\phi_{j,K},\phi_j\rangle_{\mathcal L^2(\mathcal T)}
    \bigr|,
    \qquad
    \sin\Theta_{j,K}(C)
    =
    \bigl\|
        (I-P_{j,K})\phi_j
    \bigr\|_{\mathcal L^2(\mathcal T)},
\]
where \(P_{j,K}=\phi_{j,K}\otimes\phi_{j,K}\) is the orthogonal projector onto \(\operatorname{span}\{\phi_{j,K}\}\). As shown in Section S.2.1 in the Supplementary Material \citep{supplyment}, under Assumption~\ref{ass:cov}(ii)--(iii), if \(j^e\vee j^{(\tau+1)/s} \lesssim K\), with a sufficiently small implicit constant, then uniformly over this range,
\begin{equation}
\label{equ:jiaodufanshudengjia}
    \sin\Theta_{j,K}(C)
    \asymp
    \|(I-\Pi_K)\phi_j\|_{\mathcal{L}^2(\mathcal T)}.
\end{equation}
Thus, in this range, the spectral approximation error has the same order as the eigenfunction approximation error in Assumption~\ref{ass:cov}(iii).

{More generally, the equivalence in \eqref{equ:jiaodufanshudengjia} reveals a useful insight for bounding bias in spectral estimation. Under suitable spectral regularity, the eigenfunctions of the approximated population covariance operator can achieve the same approximation error order as the underlying eigenfunctions, without an additional loss from shrinking eigengaps.  This bridges classical function smoothness/approximation conditions and spectral analysis, providing a technical route to controlling spectral approximation bias. Building on this insight, we separate population spectral approximation from stochastic perturbations. The former is controlled at the same rate as intrinsic eigenfunction approximation rate, while the latter is analyzed coefficientwise in the eigenbasis of the projected covariance operator, exploiting eigenvalue-dependent variances and individual spectral gaps. This new framework provides a key improvement over the existing analysis like \citet{zhou2024theoryfunctionalprincipalcomponent} and yields sharper spectral error bounds, suggesting a general approach to a broad class of nonparametric spectral estimation problems. }

The exponent \(e\) determines the range of spectral indices for which the sieve approximation is guaranteed.  To avoid overfitting, the sieve dimension \(K\) cannot be chosen arbitrarily large, but must be limited by the sample size. Hence, a smaller \(e\) allows larger \(j\) for a given sieve dimension \(K\).  Therefore, in the results below, \(e\) mainly determines the admissible growth range of indices \(j\), but not the convergence rate within that range.

We now state a convenient sufficient growth regime used in the minimax results below. Both the index \(j\) of an individual eigenfunction and the dimension \(J\) of a leading principal subspace are allowed to diverge with \(n\) and \(m\). For the individual eigenfunction results, it suffices that
\begin{equation}
\label{equ:elementwise-unified-frontier}
\begin{aligned}
    j=o\Bigg[
    &(nm)^{1/\{d(2s+1)+a-2\tau\}}
    \wedge (nm^2)^{1/\{d(2s+1)+2a-2\tau\}}\\
    &\wedge
    \left\{\frac{n}{\log n}\right\}^{1/(2a+2)}
    \wedge
    \left\{
        \frac{n^{s-1/2}m^{2s}}{(\log n)^{s+1/2}}
    \right\}^{1/\{2s(a+1)+2\tau+1\}}
    \Bigg].
\end{aligned}
\end{equation}
where $d = \max\{e, (\tau+1)/s\}$. For illustration, consider the second-order smoothness example with \(s=\tau=2\). Taking \(e=1\) gives \(d=3/2\), and condition~\eqref{equ:elementwise-unified-frontier} becomes $ j=o[(nm)^{2/(2a+7)} \wedge(nm^2)^{2/(4a+7)} \wedge n^{1/(2a+2)}\wedge(n^{3/2}m^4)^{1/(4a+9)}] $, up to logarithmic factors. Further specializing to $a=2$, this simplifies to  $ j=o[ n^{1/6} \wedge \left(n^{3/2}m^4\right)^{1/17} ]. $ Condition \eqref{equ:elementwise-unified-frontier} is sufficient but not necessary; see the larger range discussion in Remark~\ref{rmk:largerJ}.

For the lower bound, when \(m>j^a\), we additionally require 
\begin{equation}
\label{equ:lowerjfanwei}
\frac{j^{(2\tau-a)(a+1)}}{(nm)^{2s-a}}=o(1).
\end{equation}
This condition is needed only for certain lower-bound constructions, as specified in the results below. No additional restriction is needed when \(m\leq j^a\).  This condition is used only in the lower-bound construction of the nonparametric reconstruction error. The leading-subspace results in Section~\ref{subsec:leading-subspace} use the same conditions with \(j\) replaced by \(J\).

\subsection{Minimax theory for an individual eigenfunction}\label{sec:elementwise}
We point out the difficulties in estimating an individual eigendirection from noisy and discretely observed trajectories. The first is the local spectral identification that distinguishes neighboring eigendirections when the corresponding eigenvalues are close. The second is the  nonparametric  recovery of the target eigendirection from discrete noisy observations. These two mechanisms, together, determine the minimax rate. The following theorem formalizes this lower bound over the admissible growing-index regime.

\begin{theorem}[Lower Bound for an Individual  Eigenfunction]
\label{thm:elementwise-lower}
Let \(\mathcal Q_{\mathrm{FPCA}}\) denote the collection of all data distributions satisfying Assumptions~\ref{assump:sieve_stability} and \ref{ass:cov} with common fixed constants.  For $j$ satisfying \eqref{equ:elementwise-unified-frontier}, additionally assume \eqref{equ:lowerjfanwei} if $m>j^a$ and the second term in the rate displayed below exceeds the first. There exists a constant \(c>0\), independent of \(n,m\), and \(j\), such that
\[
\begin{aligned}
&\varliminf_{n\to\infty}
    \inf_{\widetilde\phi_j}
    \sup_{\mathcal P\in\mathcal Q_{\mathrm{FPCA}}}
    \mathbb P_{\mathcal P}\Bigg(
    \min_{\nu \in\{-1,1\}}
    \bigl\|
        \widetilde\phi_j-\nu\phi_{j,\mathcal P}
    \bigr\|^2\\
&\qquad\geq
    c\bigg[
    \frac{j^2}{n}
    \left(
        1+\frac{j^a}{m}
    \right)^2
    +
    \frac{
        j^{(2as+2\tau)/(2s+1)}
    }{
        (nm)^{2s/(2s+1)}
    }
    \left(
        1+\frac{j^a}{m}
    \right)^{2s/(2s+1)}
    \bigg]\Bigg)
    >0,
\end{aligned}
\]
where the infimum is taken over all measurable unit-norm estimators \(\widetilde\phi_j\) based on $\{ (T_{i\ell},X_{i\ell}): 1\leq i\leq n,\; 1\leq \ell\leq m_i \},$ and \(\phi_{j, \mathcal P}\) denotes a fixed oriented representative of the \(j\)th eigenfunction under \(\mathcal P\).  In the lower-bound constructions, the representatives are chosen so that the indicated pairwise inner products are nonnegative.
\end{theorem}

Specifically, under the representative choice \(a=s=\tau=2\) and for a more concise expression, the rate above becomes
\begin{equation}\label{equ:stau2f}
\frac{j^2}{n}\left(1+\frac{j^2}{m}\right)^2
+
\left\{\frac{j^{3}}{nm}
\left(1+\frac{j^2}{m}\right)\right\}^{4/5}.
\end{equation}
The two terms in the Theorem \ref{thm:elementwise-lower} arise from fundamentally different perturbation geometries. The first captures local spectral identification within a sufficiently resolved spectral neighborhood, whereas the second term concerns the recovery of the eigendirection at a resolution supported by the data. 

Specifically, the first term gives a lower bound for estimating the orientation of one principal axis of the local covariance ellipsoid. This becomes more difficult as neighboring eigenvalues get closer to each other and the ellipsoid becomes less directionally distinguishable.  The rotation construction in the lower bound quantifies the effect of eigenvalue separation. It rotates \(\phi_j\) and \(\phi_{j+1}\) by an angle \(h\). Thus, the target direction changes by order \(h\), whereas the covariance changes by order \(h(\mu_j-\mu_{j+1})\).  Under noisy discrete sampling, the effective variance levels of the two relevant score coordinates are $ \mu_j+\sigma_\varepsilon^2/{m} $ and $ \mu_{j+1}+{\sigma_\varepsilon^2}/{m}. $ Correspondingly, the KL divergence between the data distributions induced by two local orientations separated by \(h\) is bounded above by a constant multiple of $ nh^2 { (\mu_j-\mu_{j+1})^2 }/\{ (\mu_j+\sigma_\varepsilon^2/m) (\mu_{j+1}+\sigma_\varepsilon^2/m) \}. $ Hence, a statistically indistinguishable squared rotation can be chosen at the scale
\[
    h^2
    \asymp
    \frac{
        (\mu_j+\sigma_\varepsilon^2/m)
        (\mu_{j+1}+\sigma_\varepsilon^2/m)
    }{
        n(\mu_j-\mu_{j+1})^2
    }
    \asymp
    \frac{j^2}{n}
    \left(1+\frac{j^a}{m}\right)^2.
\]
Moreover, \(j^a/m\asymp m^{-1}/\mu_j\) is the inverse within-trajectory signal-to-noise ratio for the \(j\)th component. Moreover, when \(m\gtrsim j^a\), the fully observed rate \(j^2/n\) is recovered; when \(m\lesssim j^a\), sparse sampling and limited within-trajectory information inflate the error to \(j^{2a+2}/(nm^2)\).

The second term is a nonparametric reconstruction error induced by noisy discrete observations. Since one does not confine \(\phi_j\) to a common fixed finite-dimensional space, recovering it uniformly over the class requires a nonparametric trade-off. To express the geometry of the lower bound, we perturb the target eigendirection toward an \(M\)-dimensional block of unresolved higher-order eigendirections. Specifically, we consider alternatives of the form
\[
    \phi_{j,\theta}
    =
    \frac{
        \phi_j
        +
        h\sum_{M < q\leq 2M}\theta_q\phi_q
    }{
        \sqrt{1+Mh^2}
    },
    \qquad
    \theta_q\in\{-1,1\}.
\]
The resulting squared angular separation is therefore of order $Mh^2.$ Assumption~\ref{ass:cov}(iii) limits the total perturbation energy in this block, i.e. $ Mh^2 \lesssim j^{2\tau}M^{-2s}. $ When the reconstruction term dominates, the lower-bound construction bounds the KL divergence, averaged over hypercube vertices and flipped coordinates, by a constant multiple of $ nh^2 { (\mu_j-\mu_q)^2 }/\{ (\mu_j+\sigma_\varepsilon^2/m) (\mu_q+\sigma_\varepsilon^2/m) \}, $ Since \(q\asymp M\) and \(M\) is sufficiently large relative to \(j\), we have \(\mu_j-\mu_q\asymp\mu_j\). Hence, taking
\[
    h^2
    \asymp
    \frac{1}{nm\mu_j}
    \left(1+\frac{1}{m\mu_j}\right)
    \asymp
    \frac{j^a}{nm}
    \left(1+\frac{j^a}{m}\right)
\]
keeps the information bounded. Combining this information constraint with the approximation constraint gives the second term in the lower bound, characterizing the challenge of nonparametric reconstruction.

We next show that the considered estimator attains this lower bound.

\begin{theorem}[Upper Bound for an Individual Eigenfunction]
\label{thm:elementwise-expected-upper}
Under Assumptions~\ref{assump:sieve_stability} and \ref{ass:cov}, assume 
\begin{equation}
\label{equ:elementwise-expected-range}
    K \gtrsim j^e
    \vee
    j^{(\tau+1)/s},
    \qquad
    K
    =
    o\left(
        \frac{m\sqrt n}
        {j^{a+1}\sqrt{\log n}}
    \right),
    \qquad
    j^{a+1}
    \sqrt{\frac{\log n}{n}}
    =o(1),
\end{equation}
where the implicit constant in the lower bound for \(K\) is chosen sufficiently large. Then the \(j\)th empirical eigenfunction satisfies \begin{equation}
\label{equ:elementwise-expected-general}
    \|\widehat\phi_{j}-\phi_j\|_{\mathcal{L}^2(\mathcal T)}^2
    =O_p \left[
        j^{2\tau}K^{-2s}
        +
        \frac{j^2}{n}
        \left(1+\frac{j^a}{m}\right)^2 +
        \frac{Kj^a}{nm}
        \left(1+\frac{j^a}{m}\right)\right].
\end{equation}
In particular, set an integer 
$$
K\asymp (nm)^{1/(2s+1)}
        j^{(2\tau-a)/(2s+1)}
        \left(1+{j^a}/{m}\right)^{-1/(2s+1)}
$$ with suitable constant. Under the sufficient condition \eqref{equ:elementwise-unified-frontier}, the admissibility conditions in \eqref{equ:elementwise-expected-range} are naturally satisfied. Consequently,
\begin{equation}
\label{equ:elementwise-expected-oracle}
\begin{aligned}
    \|\widehat\phi_{j}-\phi_j\|_{\mathcal{L}^2(\mathcal T)}^2
    = O_p \left[
        \frac{j^2}{n}
        \left(1+\frac{j^a}{m}\right)^2 +
        \frac{
            j^{(2as+2\tau)/(2s+1)}
        }{
            (nm)^{2s/(2s+1)}
        }
        \left(1+\frac{j^a}{m}\right)^{2s/(2s+1)}\right].
\end{aligned}
\end{equation}
\end{theorem}
Combining Theorems~ \ref{thm:elementwise-lower} and \ref{thm:elementwise-expected-upper}, the proposed sieve estimator is minimax optimal for estimating \(\phi_j\) throughout the admissible growing-index regime. The optimal rate is
\begin{equation}
\label{equ:elementwise-optimal-rate}
\begin{aligned}
    \mathfrak R_j(n,m)
    \asymp\;&
    \underbrace{
        \frac{j^2}{n}
        \left(1+\frac{j^a}{m}\right)^2
    }_{\text{local spectral
identification}}
    + \ \ 
    \underbrace{
        \frac{
            j^{(2as+2\tau)/(2s+1)}
        }{
            (nm)^{2s/(2s+1)}
        }
        \left(1+\frac{j^a}{m}\right)^{2s/(2s+1)}
    }_{\text{nonparametric reconstruction error}}.
\end{aligned}
\end{equation}
Specifically, under the representative choice \(a=s=\tau=2\), the optimal rate reduces to the concise form in \eqref{equ:stau2f}. The decomposition in the upper bound mirrors the two lower-bound constructions. The first term reflects perturbations from spectral components near \(j\), whereas the second captures the cumulative contribution of spectral components outside this local neighborhood.

The closely related work of \citet{zhou2024theoryfunctionalprincipalcomponent} imposes eigenfunction derivative assumptions, which is compatible with the special case $s=2$ and $\tau=c$ of our framework. Their upper bound for random design yields a nonparametric reconstruction term of order
\[
\frac{j^{(4a+2c+2)/5}}{(nm)^{4/5}}
\left(1+\frac{j^a}{m}\right)^{4/5}.
\]
By contrast, specializing our reconstruction term to $s=2$ and $\tau=c$ gives
\[
\frac{j^{(4a+2c)/5}}{(nm)^{4/5}}
\left(1+\frac{j^a}{m}\right)^{4/5},
\]
which improves the dependence on $j$ in this term by a factor of $j^{2/5}$. Moreover, the matching lower bound established here shows that the resulting dependence on $j$ is minimax sharp over the considered model class within the admissible growing-index regime.

In the upper bound, the first term arises from the local part of the spectral perturbation expansion, involving neighboring eigendirections. Note that  $ \mathbb E \big[ \big| \langle (\widehat C-C_K)\phi_{j,K}, \phi_{q,K} \rangle \big|^2 \big] \lesssim n^{-1} (\mu_j+m^{-1})(\mu_q+m^{-1}). $ For \(q\) near \(j\), combining this coefficientwise variance bound with  \( |\mu_j-\mu_q| \gtrsim j^{-a-1}|j-q| \)  and summing the resulting normalized perturbations over the local spectral neighborhood yields the first term in the minimax rate. Thus, the upper bound accumulates the same local rotations as the lower bound. This refined perturbation analysis also explains why the fully observed contribution is \(j^2/n\), rather than the much larger rate obtained by combining a global covariance bound with a crude eigengap inequality.

The nonparametric reconstruction error in the upper bound arises from the trade-off between the first and last terms on the right-hand side of \eqref{equ:elementwise-expected-general}. At sieve resolution \(K\),  the deterministic squared approximation error satisfies  $ \sin^2\Theta_{j,K}(C) \lesssim j^{2\tau}K^{-2s}. $ The stochastic cost of increasing \(K\) is revealed by the spectral perturbation expansion.  For a sufficiently large fixed \(\kappa>1\), spectral components with \(q\geq\kappa j\) have eigengap of order $ |\mu_j-\mu_q| \asymp \mu_j \asymp j^{-a}. $ For \(q\leq j/\kappa\), the eigengap is instead of order \(\mu_q\asymp q^{-a}\); these low-index contributions are bounded separately in the spectral summation argument. For the high-index terms, normalization by \(\mu_j\) and accumulation over at most \(K\) retained sieve directions gives a \(K\)-dependent sampling contribution bounded in order by $ (K/n)/(m\mu_j) \cdot \left(1+{1}/{(m\mu_j)}\right) \asymp {Kj^a} \left(1+ {j^a}/{m}\right)/{(nm)}. $ The fully observed contribution from the same high-index range is of order \(j/n\), which is absorbed by the local spectral-identification term. Balancing the approximation and estimation errors yields the optimal rate. 

\subsection{Minimax theory for leading principal subspaces}
\label{subsec:leading-subspace}
The preceding results concern individually labelled eigendirections, for which the optimal tuning may depend on the target index. However, in  statistical modeling, such as dimension reduction, trajectory reconstruction, and subsequent analysis \citep{yao2005b,yao2005a,hall2007methodology,dou2012estimation}, constructing a separately tuned covariance estimator for each eigenfunction is impractical. Instead, a collection of principal components are typically obtained from  a common covariance estimator, and downstream analyses use the subspace spanned by such leading eigenfunctions. Estimating this principal subspace is therefore a distinct problem that carries both practical and theoretical importance.

Recall from \eqref{equ:yinruPJ} that the population and empirical projectors onto the leading \(J\)-dimensional principal subspaces are
\[
    Q_J
    =
    \sum_{j=1}^J \phi_j\otimes\phi_j,
    \qquad
    \widehat Q_J
    =
    \sum_{j=1}^J \widehat\phi_j\otimes\widehat\phi_j.
\]
Here, \(\widehat\phi_1,\ldots,\widehat\phi_J\) are obtained from the same covariance estimator and therefore share a common sieve dimension \(K\). If \(\vartheta_1,\ldots,\vartheta_J\) denote the principal angles between the population and empirical subspaces, then
\[
    \frac{1}{2J}
    \|\widehat Q_J-Q_J\|_{\mathrm{HS}}^2
    =
    \frac{1}{J}
    \sum_{j=1}^J
    \sin^2\vartheta_j
    \in[0,1].
\]

Different from the results for individually labelled eigenfunctions, estimation of the leading principal subspace has a different geometry. In particular, rotations among the first \(J\) eigendirections do not change their span and hence are invisible under the projector loss.  The sharp spectral-identification term is therefore not obtained by simply aggregating the elementwise rates over \(1\leq j\leq J\). Instead, it is determined by perturbations that move directions across the spectral boundary between the \(J\)th and \((J+1)\)st eigenvalues. Such cross-boundary rotations accumulate over many spectral coordinates, whereas rotations within either block do not contribute to the loss.

The subspace estimation problem also retains a nonparametric reconstruction error induced by discrete noisy observations. Such observations limit the resolution at which the leading eigenspace can be recovered, and eigendirections with indices near \(J\) may be simultaneously affected by perturbations along unresolved directions. Consequently, normalization by the subspace dimension does not necessarily improve the order of the reconstruction term.

We begin by identifying the information lower bound of this problem. 

\begin{theorem}[Lower Bound for Leading Principal Subspaces]
\label{thm:subspace-lower}
Let \(\mathcal Q_{\mathrm{FPCA}}\) be the collection of data distributions satisfying Assumptions~\ref{assump:sieve_stability} and \ref{ass:cov} with common fixed constants. For $J$ satisfying \eqref{equ:elementwise-unified-frontier}, additionally assume \eqref{equ:lowerjfanwei} if $m>J^a$ and the second term in the rate displayed below exceeds the first. There exists a constant \(c>0\), independent of \(n,m\), and \(J\), such that
\[
\begin{aligned}
&\varliminf_{n\to\infty}
    \inf_{\widetilde Q_J}
    \sup_{\mathcal P\in\mathcal Q_{\mathrm{FPCA}}}
    \mathbb P_{\mathcal P}\Bigg(
    \frac{1}{2J}
    \bigl\|
        \widetilde Q_J-Q_{J,\mathcal P}
    \bigr\|_{\mathrm{HS}}^2\\
&\qquad\geq
    c\bigg[
    \frac{
        J (1+\log J)
    }{n}
    \left(
        1+\frac{J^a}{m}
    \right)^2
    +
    \frac{
        J^{(2as+2\tau)/(2s+1)}
    }{
        (nm)^{2s/(2s+1)}
    }
    \left(
        1+\frac{J^a}{m}
    \right)^{2s/(2s+1)}
    \bigg]\Bigg)
    >0,
\end{aligned}
\]
where the infimum is taken over all measurable rank-\(J\) orthogonal projector estimators \(\widetilde Q_J\) based on $ \{ (T_{i\ell},X_{i\ell}): 1\leq i\leq n, 1\leq \ell\leq m_i \}, $ and \(Q_{J,\mathcal P}\) denotes the leading rank-\(J\) population spectral projector under \(\mathcal P\).
\end{theorem}

Let $m$ be large enough, it achieves the optimal rate for fully observed data \citep{wahl2022lower}.  Compared with Theorem \ref{thm:elementwise-lower}, the first term  is improved from $J^2$ to $J(1+ \log J)$.  This is because, different from a labelled eigenfunction, a principal subspace is invariant under rotations within its range. Its local identification therefore depends only on perturbations crossing the spectral boundary between the first \(J\) directions and their orthogonal complement. The lower bound formalizes this intuition by simultaneously introducing small rotations between directions \(p\in I_-\) and \(q\in I_+\), where \(L=\max\{1,\lfloor J/8\rfloor\}\), \(I_-=\{J-L+1,\ldots,J\}\), and \(I_+=\{J+1,\ldots,J+L\}\). For each pair, the lower-bound construction uses a squared rotation of order $ \rho_{pq}^2 \asymp { (\mu_p+m^{-1})(\mu_q+m^{-1}) }/{ n(\mu_p-\mu_q)^2 }. $ Hence these pairwise rotations accumulate additively under the projector loss:
\[
    \frac{1}{J}
    \sum_{p\in I_-}\sum_{q\in I_+}
    \rho_{pq}^2
    \asymp
    \frac{J(1+\log J)}{n}
    \left(
        1+\frac{J^a}{m}
    \right)^2.
\]
In particular, the $\log J$ is owing to spectral-counting, not a probability union-bound penalty.

The second term has the same order as the corresponding elementwise term in Theorem~\ref{thm:elementwise-lower}. Recall that the elementwise construction tilts a single target eigendirection toward an \(M\)-dimensional unresolved block of eigendirections. The subspace analogue uses a coupled block rotation that simultaneously tilts \(L\asymp J\) leading directions with indices near \(J\) toward the same unresolved block. The resulting squared projector separation is of order \(LMh^2\). After normalization by \(J\), it remains of order \(Mh^2\), since \(L\asymp J\). Thus, the averaging in the subspace loss is exactly offset by the number of directions simultaneously perturbed in the lower-bound construction, and the nonparametric reconstruction term retains the same order as for the eigenfunction \(\phi_J\). This differs from cross-boundary spectral identification, where rotations within the leading subspace leave the projector unchanged and therefore do not contribute to the loss.

We now establish the upper bound for the leading space estimation. A straightforward approach is to use $\frac{1}{2J} \bigl\| \widehat Q_J-Q_J \bigr\|_{\mathrm{HS}}^2 \leq \frac{1}{J} \sum_{k=1}^J \|\widehat\phi_{k}-\phi_k\|^2$. Note that an expectation analogue of \eqref{equ:elementwise-expected-general}, uniform over \(1\le k\le J\), is established in the Supplementary Material. Applying this bound with a common sieve dimension \(K\), summing over \(k\), and balancing the resulting bound in \(K\), we obtain by Markov's inequality
\begin{equation}
\label{equ:subspace-naive-upper}
    \frac{1}{2J}
    \bigl\|
        \widehat Q_J-Q_J
    \bigr\|_{\mathrm{HS}}^2
    = O_p \left[
        \frac{J^2}{n}
        \left(1+\frac{J^a}{m}\right)^2 +
        \frac{
            J^{(2as+2\tau)/(2s+1)}
        }{
            (nm)^{2s/(2s+1)}
        }
        \left(1+\frac{J^a}{m}\right)^{2s/(2s+1)}\right].
\end{equation}
It does not match the lower bound in Theorem~\ref{thm:subspace-lower}, since it ignores the fact that rotations within the leading subspace do not change its span. Moreover, practical subspace estimation uses a single covariance estimator and hence a common tuning parameter for all leading directions. Therefore, a special analysis is needed to characterize the theoretical behavior of estimating the leading principal subspace.

\begin{theorem}[Upper Bound for Leading Principal Subspaces]
\label{thm:subspace-expected-upper}
Under Assumptions~\ref{assump:sieve_stability} and \ref{ass:cov}, assume
\begin{equation}
\label{equ:subspace-expected-range}
    K
    \gtrsim
    J^e
    \vee
    J^{(\tau+1)/s},
    \qquad
    K
    =
    o\left(
        \frac{
            m\sqrt n
        }{
            J^{a+1}\sqrt{\log n}
        }
    \right),
    \qquad
    J^{a+1}
    \sqrt{\frac{\log n}{n}}
    =
    o(1),
\end{equation}
where the implicit constant in the lower bound for \(K\) is chosen sufficiently large. Then
\begin{equation}\label{equ:subspace-expected-general}
        \frac{1}{2J}
        \bigl\|
            \widehat Q_{J}-Q_J
        \bigr\|_{\mathrm{HS}}^2
    = O_p \left[
    J^{2\tau}K^{-2s}
    +
    \frac{
        J (1+\log J)
    }{n}
    \left(
        1+\frac{J^a}{m}
    \right)^2
    +
    \frac{
        KJ^a
    }{
        nm
    }
    \left(
        1+\frac{J^a}{m}
    \right)\right].
\end{equation}
In particular, set an integer
$$K
    \asymp
    (nm)^{1/(2s+1)}
    J^{(2\tau-a)/(2s+1)}
    \left(
        1+{J^a}/{m}
    \right)^{-1/(2s+1)},$$
with suitable constant. Under the sufficient condition \eqref{equ:elementwise-unified-frontier}, with \(j\) replaced by \(J\), the admissibility conditions in \eqref{equ:subspace-expected-range} are naturally satisfied. Consequently,
$$
        \frac{1}{2J}
        \bigl\|
            \widehat Q_{J}-Q_J
        \bigr\|_{\mathrm{HS}}^2
     = 
    O_p \left[
    \frac{
        J (1+\log J)
    }{n}
    \left(
        1+\frac{J^a}{m}
    \right)^2
    +
    \frac{
        J^{(2as+2\tau)/(2s+1)}
    }{
        (nm)^{2s/(2s+1)}
    }
    \left(
        1+\frac{J^a}{m}
    \right)^{2s/(2s+1)} \right].
$$
\end{theorem}

Combining Theorems~\ref{thm:subspace-expected-upper} and \ref{thm:subspace-lower}, the minimax-optimal rate for estimating the rank-\(J\) leading principal subspace is
$$
\begin{aligned}
    \mathfrak R_J^{\mathrm{space}}(n,m)
    \asymp
    \underbrace{
        \frac{
            J (1+\log J) 
        }{n}
        \left(
            1+\frac{J^a}{m}
        \right)^2
    }_{\text{cross-boundary spectral identification}}
    +\ 
    \underbrace{
        \frac{
            J^{(2as+2\tau)/(2s+1)}
        }{
            (nm)^{2s/(2s+1)}
        }
        \left(
            1+\frac{J^a}{m}
        \right)^{2s/(2s+1)}
    }_{\text{nonparametric reconstruction error}}.
\end{aligned}
$$

As discussed above, the first term represents the error of distinguishing rotations that cross the spectral boundary between the first \(J\) eigendirections and their orthogonal complement. This distinction is explicit in the spectral projector expansion. Let \(Q_{J,K}\) be the leading rank-\(J\) projector of \(C_K\). Once the empirical leading cluster is correctly localized,
\begin{equation}
\label{equ:subspace-cross-boundary-expansion}
    \frac{1}{2J}
    \bigl\|
        \widehat Q_J-Q_{J,K}
    \bigr\|_{\mathrm{HS}}^2
    \lesssim
    \frac{1}{J}
    \sum_{p\leq J}
    \sum_{J<q\leq K}
    \frac{
        \bigl|
            \langle
                (\widehat C-C_K)\phi_{p,K},
                \phi_{q,K}
            \rangle
        \bigr|^2
    }{
        (\mu_{p,K}-\mu_{q,K})^2
    },
\end{equation}
where \((\mu_{p,K},\phi_{p,K})\) are the eigenpairs of \(C_K\). An important point is that only pairs satisfying \(p\leq J<q\) appear. Intra-block perturbations within either the leading subspace or its orthogonal complement merely rotate the corresponding basis and leave the spectral projector unchanged. The projector error therefore captures only inter-block perturbations connecting the two spectral groups.

The second term represents nonparametric reconstruction and has the same order as the corresponding elementwise reconstruction error for \(\phi_J\). Since the normalized projector loss averages the squared principal angles over the first \(J\) directions, one might expect the subspace reconstruction error to be smaller than the corresponding error for \(\phi_J\). Such an improvement would occur if the reconstruction difficulty were concentrated in only a few leading directions. However, all indices in the upper portion of the leading subspace, such as \(J/2\leq k\leq J\), have approximation and discretization difficulties of the same order as those for \(k=J\). Thus, order \(J\) leading eigendirections can be simultaneously affected by perturbations along unresolved directions. Their cumulative contribution offsets the factor \(1/J\) in the normalized projector loss.

\begin{remark}\label{rmk:largerJ}
Condition~\eqref{equ:elementwise-unified-frontier} is a convenient sufficient condition ensuring that the balancing choice of \(K\) satisfies the admissibility conditions in \eqref{equ:elementwise-expected-range} and \eqref{equ:subspace-expected-range}. For example, when \(a=s=\tau=2\) and \(e=1\), a refined analysis shows that the final rates in Theorems~\ref{thm:subspace-lower} and~\ref{thm:subspace-expected-upper} also hold for a larger range $ J=o\{n^{1/4}\wedge(nm)^{1/6}\wedge(nm^2)^{2/15}\} $ with balancing $K$. The same statement holds for Theorems~\ref{thm:elementwise-lower} and~\ref{thm:elementwise-expected-upper} with \(J\) replaced by \(j\). Moreover, the balancing choice of \(K\) is not necessary when the nonparametric reconstruction term is no longer rate-determining. For example, suppose that $ e=1$, $s>\tau+a+\frac12, $ and we can consider $J$ with
\begin{equation}\label{equ:largerJ}
 J=o\left[ \left(\frac{n}{\log n}\right)^{\frac1{2a+2}} \wedge \left(\frac{nm^2}{\log n}\right)^{\frac1{2a+4}} \right]. 
\end{equation}
Under these conditions, Theorem~\ref{thm:subspace-expected-upper} requires \(K\gtrsim J\), with a sufficiently large implicit constant. If the balancing choice of \(K\) is only $o(J)$, one may take \(K\asymp J\).  Under the displayed condition, this choice remains admissible. In this case, increasing \(K\) to order \(J\) reduces the approximation error, while the stochastic term associated with \(K\asymp J\) is dominated by the spectral-identification term. Hence the same upper bound continues to hold under the larger range~\eqref{equ:largerJ}. The corresponding lower bound also holds under this larger range. See the Supplementary Material \citep{supplyment} for detailed proof. The same argument applies to the estimation of a labelled eigenfunction, with \(j\) in place of \(J\).
\end{remark}

\subsection{Phase transitions}
\label{subsec:sampling-regimes}
The minimax rates derived in the preceding subsections characterize how the estimation risks for eigenfunctions and principal subspaces depend jointly on the sample size, sampling frequency, and target spectral order. This subsection gives explicit rate expressions for the sparse-to-dense sampling regimes. The sampling intervals below are considered over their intersections with the admissible range, which includes the extension in Remark~\ref{rmk:largerJ}. We first consider the estimation of an individually labelled eigenfunction.

\begin{corollary}[Convergence Regimes for an Individual Eigenfunction]
\label{cor:elementwise-sampling-regimes}
Suppose that the conditions of Theorems~\ref{thm:elementwise-lower} and \ref{thm:elementwise-expected-upper} hold or the conditions in Remark \ref{rmk:largerJ} hold.

\begin{enumerate}[(i)]
\item If $ j^{2s+1-\tau} = o\bigl(n^{1/2}\bigr), $ the rate admits the following four-regime characterization over the admissible sampling range
$$
    \mathfrak R_j(n,m)
    \asymp
    \begin{cases}
        \displaystyle
        \frac{j^{2a+2}}{nm^2},
        &
        m
        \lesssim
        n^{-1/2}j^{a+2s+1-\tau},
        \\[12pt]
        \displaystyle
        \frac{
            j^{(4as+2\tau)/(2s+1)}
        }{
            (nm^2)^{2s/(2s+1)}
        },
        &
        n^{-1/2}j^{a+2s+1-\tau}
        \lesssim
        m
        \lesssim
        j^a,
        \\[14pt]
        \displaystyle
        \frac{
            j^{(2as+2\tau)/(2s+1)}
        }{
            (nm)^{2s/(2s+1)}
        },
        &
        j^a
        \lesssim
        m
        \lesssim
        n^{1/(2s)}
        j^{a-2+(\tau-1)/s},
        \\[14pt]
        \displaystyle
        \frac{j^2}{n},
        &
        m
        \gtrsim
        n^{1/(2s)}
        j^{a-2+(\tau-1)/s}.
    \end{cases}
$$

\item  If \(j^{2s+1-\tau}\gtrsim n^{1/2}\), there are two convergence regimes
$$
    \mathfrak R_j(n,m)
    \asymp
    \begin{cases}
        \displaystyle
        \frac{j^{2a+2}}{nm^2},
        &
        m\lesssim j^a,
        \\[11pt]
        \displaystyle
        \frac{j^2}{n},
        &
        m\gtrsim j^a.
    \end{cases}
$$
\end{enumerate}
\end{corollary}

Corollary~\ref{cor:elementwise-sampling-regimes} shows that the number of convergence regimes is determined by the relative magnitude of the spectral-identification and nonparametric reconstruction errors near the sampling threshold \(m^{-1}\asymp\mu_j\), or equivalently \(m\asymp j^a\), and Figure~\ref{fig:pointwise-phase-diagram} graphically illustrates such transitions.

At \(m\asymp j^a\), the two components are of respective orders ${j^2}/{n} \text{ and }{ j^{2\tau/(2s+1)} }/{ n^{2s/(2s+1)} }.$ In the case that $ j^{2s+1-\tau}=o(n^{1/2}), $ i.e. ${j^2}/{n} = o(j^{2\tau/(2s+1)} /{ n^{2s/(2s+1)} })$,  nonparametric reconstruction error is rate-determining in a neighborhood of \(m\asymp j^a\). Away from this boundary and for sufficiently sparse sampling, the rate is dominated by local spectral identification error. When \(m\gtrsim j^a\), the spectral-identification term approaches the fully observed rate \(j^2/n\), whereas the reconstruction error decreases with \(m\) and becomes very small. Spectral identification therefore becomes rate-determining again under sufficiently dense sampling.  Consequently, whenever nonparametric reconstruction dominates near \(m\asymp j^a\), the two error terms cross once on each side of \(m\asymp j^a\),  giving rise to four asymptotically separated convergence regimes when the transition scales lie in the admissible sampling range. This phase diagram also connects several existing theories.  Specifically, in this phase diagram, the theory of \citet{hall2006properties} with fixed $m$ and $j$ corresponds to the second/third regime, the fully observed growing-index theory of \citet{hall2007methodology} corresponds to the fourth regime, while the fixed-index theory of \citet{cai2010nonparametric} across arbitrary sampling frequencies spans the third and fourth regimes. The first regime is new and arises specifically from the diverging spectral index in sparse sampling.

These four regimes admit a direct interpretation as the sampling frequency increases from sparse to dense. Under extremely sparse sampling, the main difficulty is identifying the orientation of the \(j\)th principal axis because the neighboring principal components are only weakly resolved within each trajectory. As \(m\) increases, the local spectral-identification term becomes smaller than the nonparametric reconstruction term, while the target score remains insufficiently resolved. After crossing the sampling threshold \(m\asymp j^a\), the target score becomes sufficiently resolved within each trajectory, but nonparametric reconstruction of the eigendirection remains rate-determining.  Finally, as the sampling becomes sufficiently dense, the nonparametric reconstruction error becomes negligible, and the spectral-identification term becomes dominant again. The fully observed rate $j^2/n$ is then attained, and further increases in \(m\) yield no rate improvement.

In the case that $ j^{2s+1-\tau}\gtrsim n^{1/2}, $ the local spectral-identification error is already no smaller than the nonparametric reconstruction error at \(m\asymp j^a\), where the latter is most competitive.  Hence, local spectral identification remains rate-determining throughout the entire sampling range, and the two intermediate reconstruction regimes disappear.  The only remaining transition occurs at \(m\asymp j^a\), separating the regimes in which the target score is poorly and sufficiently resolved within each trajectory.  To let \(j^{2s+1-\tau}\gtrsim n^{1/2}\) can occur as \(j\) increases, we require \(\tau<2s+1\). Under Assumption~\ref{ass:cov}(iii), this means that the sieve approximation error for higher-order eigenfunctions does not increase too rapidly with the spectral index. Consequently, at sufficiently large \(j\), the shrinking eigengaps make local spectral identification more difficult than nonparametric reconstruction throughout the sampling range, leaving only the two spectral-identification regimes.

The following result is for the leading principal subspace.

\begin{corollary}[Convergence Regimes for a Leading Principal Subspace]
\label{cor:subspace-sampling-regimes}
Suppose that the conditions of Theorems~\ref{thm:subspace-lower} and \ref{thm:subspace-expected-upper} hold  or the conditions in Remark \ref{rmk:largerJ} hold.

\begin{enumerate}[(i)]
\item If $J^{s+1/2-\tau} (1+\log J)^{s+1/2} = o\bigl(n^{1/2}\bigr),$ the rate admits the following four-regime characterization over the admissible sampling range
$$
    \mathfrak R_J^{\mathrm{space}}(n,m)
    \asymp
    \begin{cases}
        \displaystyle
        \frac{
            J^{2a+1}(1+\log J)
        }{
            nm^2
        },
        &
        m
        \lesssim
        n^{-1/2}
        J^{a+s+1/2-\tau}
        (1+\log J)^{(2s+1)/2},
        \\[14pt]
        \displaystyle
        \frac{
            J^{(4as+2\tau)/(2s+1)}
        }{
            (nm^2)^{2s/(2s+1)}
        },
        &
        n^{-1/2}
        J^{a+s+1/2-\tau}
        (1+\log J)^{(2s+1)/2}
        \lesssim
        m
        \lesssim
        J^a,
        \\[15pt]
        \displaystyle
        \frac{
            J^{(2as+2\tau)/(2s+1)}
        }{
            (nm)^{2s/(2s+1)}
        },
        &
        J^a
        \lesssim
        m
        \lesssim
        n^{1/(2s)}
        J^{a-1+(\tau-1/2)/s}
        (1+\log J)^{-(2s+1)/(2s)},
        \\[15pt]
        \displaystyle
        \frac{
            J(1+\log J)
        }{
            n
        },
        &
        m
        \gtrsim
        n^{1/(2s)}
        J^{a-1+(\tau-1/2)/s}
        (1+\log J)^{-(2s+1)/(2s)}.
    \end{cases}
$$

\item  If $J^{s+1/2-\tau} (1+\log J)^{s+1/2} \gtrsim n^{1/2},$ there are two convergence regimes
$$
    \mathfrak R_J^{\mathrm{space}}(n,m)
    \asymp
    \begin{cases}
        \displaystyle
        \frac{
            J^{2a+1}(1+\log J)
        }{
            nm^2
        },
        &
        m\lesssim J^a,
        \\[13pt]
        \displaystyle
        \frac{
            J(1+\log J)
        }{
            n
        },
        &
        m\gtrsim J^a.
    \end{cases}
$$
\end{enumerate}
\end{corollary}

The leading-subspace problem exhibits the same qualitative phase-transition structure, but the four-regime configuration persists to a larger spectral order, which is also illustrated in Figure~\ref{fig:subspace-phase-diagram}  At the boundary \(m\asymp J^a\), the spectral-identification error is ${J(1+\log J)}/{n},$ rather than \(J^2/n\) for the labelled eigenfunction \(\phi_J\).  As discussed above, this reduction is a consequence of the projector geometry: rotations within the leading \(J\)-dimensional subspace leave the projector unchanged, hence only perturbations crossing the spectral boundary contribute to the loss. Consequently, the four-regime condition $ J^{s+1/2-\tau} (1+\log J)^{s+1/2} = o(n^{1/2}) $ is substantially weaker than its elementwise counterpart.  The admissible range considered here also includes the extension in Remark~\ref{rmk:largerJ}. There can be a broad range of indices for which the individully labelled eigenfunction has already entered the two-regime phase, while the leading principal subspace still exhibits all four regimes.

\section{Application to Functional Linear Regression}\label{sec:functional-linear-regression}

Functional principal component analysis is widely used for dimension reduction and as a regularization tool in functional regression. Functional linear regression provides a canonical example \citep{yao2005b,hall2007methodology,dou2012estimation}. In this section, we analyze FPCA-based functional linear regression with noisy and discretely observed predictors.

Consider the functional linear model
\begin{equation}
    Y_i
    =
    \int_{\mathcal T} X_i(t)\beta(t) dt
    +
    \zeta_i,
    \qquad
    1\leq i\leq n,
    \label{equ:flr-model}
\end{equation}
where the trajectories are discretely observed according to \eqref{equ:lisanx}. The response errors \(\zeta_i\) are i.i.d., centered, have bounded variance, and are independent of all other random variables.

\begin{assumption}[Slope Regularity]
\label{ass:flr-slope}
The slope function satisfies $ \beta = \sum_{k=1}^{\infty}\beta_k\phi_k, $ with $ |\beta_k| \lesssim k^{-b}, $ uniformly over $k\geq1,$ where $b>\max\{\tau+1,e+1/2, a/2+1\}$.
\end{assumption}

The exponent \(b\) characterizes the regularity of the slope relative to the covariance eigenbasis. A larger value of \(b\) corresponds to faster decay of the coefficients along higher-order eigendirections and hence to less weight on these directions, leading to an easier estimation and prediction problem \citep{hall2007methodology,dou2012estimation,zhou2023functional}.

To evaluate an estimator \(\widetilde\beta\), we consider both the squared \(\mathcal L^2\)-loss and the prediction loss \citep{hall2007methodology,cai2012minimax}. The former is $ \|\widetilde\beta-\beta\|_{\mathcal L^2(\mathcal T)}^2, $ which directly measures error in recovering the slope function. The prediction loss is $ \|C^{1/2}(\widetilde\beta-\beta)\|_{\mathcal L^2(\mathcal T)}^2 = \mathbb E [ \{ \int_{\mathcal T} X(t) \bigl(\widetilde\beta(t)-\beta(t)\bigr)\,dt \}^2 ], $ where the expectation is taken with respect to an independent copy \(X\). Thus, the latter is the excess mean-squared prediction error for a new trajectory.

The minimax theory for functional linear regression with discretely observed predictors has been a longstanding open problem and was recently established by \citet{yanflr} for a pooling-ridge estimation within RKHS framework. Under \(\mu_k\asymp k^{-a}\) and \(|\beta_k|\lesssim k^{-b}\), the minimax rates under the squared \(\mathcal L^2\)-loss and prediction loss are, respectively,
\begin{equation}
\label{equ:flrlowerl2}
    n^{-(2b-1)/(a+2b)}
    +
    (nm)^{-(2b-1)/(2a+2b)},
    \
    n^{-(a+2b-1)/(a+2b)}
    +
    (nm)^{-(a+2b-1)/(2a+2b)}.
\end{equation}
Relatedly, \citet{cai2026costdiscretizationfunctionallinear} derive the minimax rates for discretely observed functional linear regression under trigonometric eigenfunctions, and develop adaptive procedures for unknown eigenvalues.  This isolates the effect of discretization and provides useful insight into the resulting sampling phase transitions, while setting aside the challenge of estimating the eigenbasis from the discretely observed data. It therefore remains to determine whether the classical FPCA-based FLR procedure, with empirical eigenfunctions estimated from data, can attain the corresponding minimax rates.

We now introduce the considered FPCA-based estimator for functional linear regression.  Define the cross-covariance function between the response and the functional predictor by
\[
    h(t)
    :=
    \mathbb E[YX(t)]
    =
    \int_{\mathcal T}
        C(t,u)\beta(u)\,du.
\]
Equivalently, \(h=C\beta\). For a sieve dimension \(K\), we estimate \(h\) directly from the discrete observations by
$$
\widehat h =
    \frac{1}{n}
    \sum_{i=1}^n
    \frac1{m_i}\sum_{j=1}^{m_i}
    Y_iX_{ij}g_{K,T_{ij}}.
$$
Obviously, we have \( \mathbb E\widehat h = \Pi_KC\beta. \) {Recall \(\widehat C\) as the covariance estimator in \eqref{equ:covariance_estimator}, and let \( \widehat Q_J \) be the spectral projector onto its first \(J\) empirical eigendirections. We define the same-sample empirical principal-component estimator by
\begin{equation}
    \widehat\beta
    =
    \bigl(
        \widehat Q_J
        \widehat C
        \widehat Q_J
    \bigr)^\dagger
    \widehat Q_J\widehat h=
    \sum_{k=1}^J
    \frac{
        \langle
            \widehat h,
            \widehat\phi_k
        \rangle_{\mathcal{L}^2(\mathcal T)}
    }{
        \widehat\mu_k
    }
    \widehat\phi_k, 
    \label{equ:flr-estimator}
\end{equation}
where the last equality holds when the empirical eigenvalues are nonzero and simple.  }

The following theorem establishes the convergence rates of the FPCA-based estimator under both the squared \(\mathcal{L}^2\)-loss and the prediction loss. Note that Assumption~\ref{ass:cov}(iii) is monotone in $e$: if it holds for some $e\ge1$, it also holds for any larger $e$. Thus, if needed, a larger $e$ can be adopted in the following without strengthening the approximation assumption.

\begin{theorem}[Upper Bound for Functional Linear Regression]
\label{thm:flr-upper}
Under Assumptions~\ref{assump:sieve_stability}, \ref{ass:cov}, and \ref{ass:flr-slope}, suppose that $s>b/e-1/(2e)$, $es\geq\tau+1$, and $2a+1>e$. Let \(J\to\infty\) satisfying $J^{2a+2}=o\left({n}/{\log n}\right).$ Let \(K\) be a deterministic integer satisfying
\begin{equation}
 \label{equ:flr-K-balance}
J^e \lesssim K
\lesssim \{
mJ^{2b-1}
\wedge
\left(
    J^{2a+1}+mJ^{a+1}
\right) \}/ (1+\log J)^2,
\end{equation}
for suitable implicit constants; when \(es=\tau+1\), the constant in \(K\gtrsim J^e\) is sufficiently large. The factor \((1+\log J)^2\) may be omitted when \(b\neq a+1\). We also assume
\begin{equation}
\label{equ:flr-K-upper-range}
K=o\left[
\left\{
\frac{
nm^2
}{
J^{2a+1}\log n
}
\right\}^{1/2}
\wedge
\frac{nm}{J^{2a+2}\log n}
\wedge
\frac{n}{J^{a+1}\log^2 n}
\right].
\end{equation}
Then, the estimator \(\widehat\beta\) in \eqref{equ:flr-estimator} satisfies
\begin{equation}
        \big\|
        \widehat\beta-\beta
    \big\|_{\mathcal{L}^2(\mathcal T)}^2
    =
    O_p\left[
        J^{1-2b}
        +
        \frac{J^{a+1}}{n}
        +
        \frac{J^{2a+1}}{nm}
    \right],\label{equ:flr-l2-J-rate}
\end{equation}
and
\begin{equation}
    \big\|
        C^{1/2}
        (\widehat\beta-\beta)
    \big\|_{\mathcal{L}^2(\mathcal T)}^2
    =
    O_p\left[
        J^{1-a-2b}
        +
        \frac{J}{n}
        +
        \frac{J^{a+1}}{nm}
    \right].
    \label{equ:flr-prediction-J-rate}
\end{equation}
In particular, set the positive integer
$$J \asymp
        n^{1/(a+2b)}
        \wedge
        (nm)^{1/(2a+2b)},$$
with suitable constant. Then the admissible range of \(K\) satisfying conditions~\eqref{equ:flr-K-balance}--\eqref{equ:flr-K-upper-range} is nonempty, and hence the estimator \(\widehat\beta\) satisfies
$$
    \big\|
        \widehat\beta-\beta
    \big\|_{\mathcal{L}^2(\mathcal T)}^2
    =
    O_p\Big[
        n^{-(2b-1)/(a+2b)}
        +
        (nm)^{-(2b-1)/(2a+2b)}
    \Big],
$$
and
$$
    \big\|
        C^{1/2}
        (\widehat\beta-\beta)
    \big\|_{\mathcal{L}^2(\mathcal T)}^2
    =
    O_p\Big[
        n^{-(a+2b-1)/(a+2b)}
        +
        (nm)^{-(a+2b-1)/(2a+2b)}
    \Big].
$$
\end{theorem}

The three terms in \eqref{equ:flr-prediction-J-rate} have distinct origins. The term \(J^{1-a-2b}\) is the truncation bias of the slope function. The term \(J/n\) is the subject-level sampling error and coincides with the contribution under fully observed predictors, whereas \(J^{a+1}/(nm)\) is the additional error induced by noisy discrete sampling. Under the \(\mathcal L^2\)-loss, removing the covariance weighting amplifies errors along the weaker eigendirections by a factor of order \(\mu_J^{-1}\asymp J^a\). Consequently, the truncation bias increases to \(J^{1-2b}\), the subject-level sampling error and the discrete-sampling error become \(J^{a+1}/n\) and \(J^{2a+1}/(nm)\). Balancing the truncation and stochastic errors yields the final rates. These rates coincide with the minimax lower bounds in \eqref{equ:flrlowerl2} established by \citet{yanflr}, thereby proving the minimax optimality of the FPCA-based estimator under both losses.

The bounds in Theorem~\ref{thm:flr-upper} exhibit the same phase transition at $m\asymp n^{a/(a+2b)}$ in terms of both the squared \(\mathcal L^2\)-loss and the prediction loss. When $ m\gtrsim n^{a/(a+2b)}, $ the regression-relevant principal component scores are sufficiently resolved within each trajectory, and the rates reduce to $ \|\widehat\beta-\beta\|_{\mathcal L^2(\mathcal T)}^2 = O_p\left( n^{-(2b-1)/(a+2b)} \right), $ and $ \|C^{1/2}(\widehat\beta-\beta)\|_{\mathcal L^2(\mathcal T)}^2 = O_p\left( n^{-(a+2b-1)/(a+2b)} \right). $ These coincide with the fully observed benchmarks \citep{hall2007methodology,dou2012estimation}. Thus, in the dense regime, the number of independent trajectories becomes the statistical bottleneck, and further increasing the sampling frequency yields no rate improvement. When $ m\lesssim n^{a/(a+2b)}, $ the within-trajectory information is insufficient for regression at the fully observed resolution, and the error induced by discrete sampling becomes rate-determining, hence $ \|\widehat\beta-\beta\|_{\mathcal L^2(\mathcal T)}^2 = O_p\left( (nm)^{-(2b-1)/(2a+2b)} \right), $ and $ \|C^{1/2}(\widehat\beta-\beta)\|_{\mathcal L^2(\mathcal T)}^2 = O_p\left( (nm)^{-(a+2b-1)/(2a+2b)} \right). $ Hence, in the sparse regime, the rates depend explicitly on the total number \(nm\) of predictor measurements.

The transition boundary also highlights the distinct roles of the spectral decay \(a\) and the slope regularity \(b\). A larger \(b\) implies faster decay of the slope coefficients along higher-order eigendirections, thereby reducing the effective complexity of the regression problem. This leads to faster convergence and lowers the sampling frequency required to attain the fully observed benchmark. A larger \(a\), by contrast, corresponds to faster decay of the principal component variances. Consequently, higher-order components carry weaker within-trajectory signals and require more observations to be reliably resolved. The sampling frequency required to attain the fully observed rate therefore increases with \(a\), shifting the phase-transition boundary toward denser sampling. It gives rise to a curse-of-smoothness phenomenon for functional inverse problems. Faster spectral decay reduces the contribution of higher-order predictor directions, but at the same time makes these directions more difficult to recover from discrete observations.  Its effect on the convergence rate depends on the loss. Under the \(\mathcal L^2\)-loss, a larger \(a\) increases the degree of ill-posedness and slows slope estimation. Under the prediction loss, however, errors along weak predictor directions are downweighted by the covariance operator, so faster spectral decay can improve the fully observed rate. In both cases, a larger \(a\) increases the sampling frequency required to reach the dense regime.

\section{Conclusion and Discussion}
This work establishes the first minimax theory for functional principal component analysis from discretely observed data that accommodates both individually labelled eigendirections with diverging indices and leading principal subspaces of increasing dimensions, which has been open for nearly two decades. The impact of discrete sampling is determined jointly by the number of trajectories, the within-trajectory sampling frequency, and the spectral order of the target. The optimal rates are governed by spectral identification and nonparametric reconstruction, whose interplay gives rise to either four or two distinct sampling regimes across the range from sparse to dense designs.

{Our analysis develops a general framework for sharp minimax theory in nonparametric spectral estimation, including lower bound constructions, and refined analysis of nonparametric approximation and stochastic spectral perturbations. For the lower bounds, we construct rotations between neighboring eigendirections and simultaneous perturbations across the spectral boundary of leading subspaces for the spectral identification errors. Complementary constructions based on orthogonal rotations toward high-frequency blocks, using Assouad hypercubes and sign-matrix packings, characterize nonparametric reconstruction errors. For the upper bounds, we adopt a new approach that separates population spectral approximation from stochastic perturbations. The former is bounded by the same rate as the intrinsic eigenfunction approximation rate, while the latter is analyzed by exploiting eigenvalue-dependent variances and individual spectral gaps in the eigenbasis of the projected covariance.  Broadly speaking, the new techniques may facilitate sharp minimax analyses of other FPCA estimation methods, including local polynomial \citep{zhou2023functional} and RKHS-based \citep{cai2010nonparametric} approaches, which deserves further investigation. }

An important further direction is to understand how this spectral analysis informs downstream FPCA-based statistical procedures. Our functional linear regression result provides one example: the classical FPCA-based estimator attains the optimal rate when the principal components are estimated from discretely observed data. This suggests that sharper spectral analysis can lead to better theoretical guarantees for downstream procedures. It would therefore be interesting to investigate whether analogous results can be established for various FPCA-based generalized/single-index/nonlinear type of regression, recovery, classification, and so on.

\bibliographystyle{imsart-nameyear}
\bibliography{bibliography}

\begin{thebibliography}{43}

\bibitem[\protect\citeauthoryear{Amini and
  Wainwright}{2012}]{10.1214/12-AOS1033}
\begin{barticle}[author]
\bauthor{\bsnm{Amini},~\bfnm{Arash~A.}\binits{A.~A.}} \AND
  \bauthor{\bsnm{Wainwright},~\bfnm{Martin~J.}\binits{M.~J.}}
(\byear{2012}).
\btitle{{Sampled forms of functional PCA in reproducing kernel Hilbert
  spaces}}.
\bjournal{The Annals of Statistics}
\bvolume{40}
\bpages{2483 -- 2510}.
\end{barticle}
\endbibitem

\bibitem[\protect\citeauthoryear{Aneiros et~al.}{2019}]{aneiros2019recent}
\begin{barticle}[author]
\bauthor{\bsnm{Aneiros},~\bfnm{Germ{\'a}n}\binits{G.}},
  \bauthor{\bsnm{Cao},~\bfnm{Ricardo}\binits{R.}},
  \bauthor{\bsnm{Fraiman},~\bfnm{Ricardo}\binits{R.}},
  \bauthor{\bsnm{Genest},~\bfnm{Christian}\binits{C.}} \AND
  \bauthor{\bsnm{Vieu},~\bfnm{Philippe}\binits{P.}}
(\byear{2019}).
\btitle{Recent advances in functional data analysis and high-dimensional
  statistics}.
\bjournal{Journal of Multivariate Analysis}
\bvolume{170}
\bpages{3--9}.
\end{barticle}
\endbibitem

\bibitem[\protect\citeauthoryear{Belhakem et~al.}{2025}]{belhakem2025minimax}
\begin{barticle}[author]
\bauthor{\bsnm{Belhakem},~\bfnm{Ryad}\binits{R.}},
  \bauthor{\bsnm{Picard},~\bfnm{Franck}\binits{F.}},
  \bauthor{\bsnm{Rivoirard},~\bfnm{Vincent}\binits{V.}} \AND
  \bauthor{\bsnm{Roche},~\bfnm{Angelina}\binits{A.}}
(\byear{2025}).
\btitle{Minimax estimation of functional principal components from noisy
  discretized functional data}.
\bjournal{Scandinavian Journal of Statistics}
\bvolume{52}
\bpages{38-80}.
\end{barticle}
\endbibitem

\bibitem[\protect\citeauthoryear{Brenner and
  Scott}{2008}]{brenner2008mathematical}
\begin{bbook}[author]
\bauthor{\bsnm{Brenner},~\bfnm{Susanne~C}\binits{S.~C.}} \AND
  \bauthor{\bsnm{Scott},~\bfnm{L~Ridgway}\binits{L.~R.}}
(\byear{2008}).
\btitle{The mathematical theory of finite element methods}.
\bpublisher{Springer}.
\end{bbook}
\endbibitem

\bibitem[\protect\citeauthoryear{Cai and Hall}{2006}]{cai2006prediction}
\begin{barticle}[author]
\bauthor{\bsnm{Cai},~\bfnm{T.~Tony}\binits{T.~T.}} \AND
  \bauthor{\bsnm{Hall},~\bfnm{Peter}\binits{P.}}
(\byear{2006}).
\btitle{{Prediction in functional linear regression}}.
\bjournal{The Annals of Statistics}
\bvolume{34}
\bpages{2159 -- 2179}.
\end{barticle}
\endbibitem

\bibitem[\protect\citeauthoryear{Cai and
  Li}{2026}]{cai2026costdiscretizationfunctionallinear}
\begin{barticle}[author]
\bauthor{\bsnm{Cai},~\bfnm{T~Tony}\binits{T.~T.}} \AND
  \bauthor{\bsnm{Li},~\bfnm{Yicheng}\binits{Y.}}
(\byear{2026}).
\btitle{The Cost of Discretization in Functional Linear Regression: Minimax
  Rates and Adaptation}.
\bjournal{arXiv preprint arXiv:2607.09350}.
\end{barticle}
\endbibitem

\bibitem[\protect\citeauthoryear{Cai and Yuan}{2010}]{cai2010nonparametric}
\begin{barticle}[author]
\bauthor{\bsnm{Cai},~\bfnm{Tony}\binits{T.}} \AND
  \bauthor{\bsnm{Yuan},~\bfnm{Ming}\binits{M.}}
(\byear{2010}).
\btitle{Nonparametric covariance function estimation for functional and
  longitudinal data}.
\bjournal{University of Pennsylvania and Georgia inistitute of technology}.
\end{barticle}
\endbibitem

\bibitem[\protect\citeauthoryear{Cai and Yuan}{2011}]{cai2011optimal}
\begin{barticle}[author]
\bauthor{\bsnm{Cai},~\bfnm{T~Tony}\binits{T.~T.}} \AND
  \bauthor{\bsnm{Yuan},~\bfnm{Ming}\binits{M.}}
(\byear{2011}).
\btitle{Optimal estimation of the mean function based on discretely sampled
  functional data: Phase transition}.
\bjournal{The Annals of Statistics}
\bvolume{39}
\bpages{2330--2355}.
\end{barticle}
\endbibitem

\bibitem[\protect\citeauthoryear{Cai and Yuan}{2012}]{cai2012minimax}
\begin{barticle}[author]
\bauthor{\bsnm{Cai},~\bfnm{T~Tony}\binits{T.~T.}} \AND
  \bauthor{\bsnm{Yuan},~\bfnm{Ming}\binits{M.}}
(\byear{2012}).
\btitle{Minimax and adaptive prediction for functional linear regression}.
\bjournal{Journal of the American Statistical Association}
\bvolume{107}
\bpages{1201--1216}.
\end{barticle}
\endbibitem

\bibitem[\protect\citeauthoryear{Chen and Lei}{2015}]{chen2015localized}
\begin{barticle}[author]
\bauthor{\bsnm{Chen},~\bfnm{Kehui}\binits{K.}} \AND
  \bauthor{\bsnm{Lei},~\bfnm{Jing}\binits{J.}}
(\byear{2015}).
\btitle{Localized functional principal component analysis}.
\bjournal{Journal of the American Statistical Association}
\bvolume{110}
\bpages{1266--1275}.
\end{barticle}
\endbibitem

\bibitem[\protect\citeauthoryear{Delaigle, Hall and
  Bathia}{2012}]{delaigle2012componentwise}
\begin{barticle}[author]
\bauthor{\bsnm{Delaigle},~\bfnm{A}\binits{A.}},
  \bauthor{\bsnm{Hall},~\bfnm{P}\binits{P.}} \AND
  \bauthor{\bsnm{Bathia},~\bfnm{N}\binits{N.}}
(\byear{2012}).
\btitle{Componentwise classification and clustering of functional data}.
\bjournal{Biometrika}
\bvolume{99}
\bpages{299--313}.
\end{barticle}
\endbibitem

\bibitem[\protect\citeauthoryear{Dou, Pollard and
  Zhou}{2012}]{dou2012estimation}
\begin{barticle}[author]
\bauthor{\bsnm{Dou},~\bfnm{Winston~Wei}\binits{W.~W.}},
  \bauthor{\bsnm{Pollard},~\bfnm{David}\binits{D.}} \AND
  \bauthor{\bsnm{Zhou},~\bfnm{Harrison~H}\binits{H.~H.}}
(\byear{2012}).
\btitle{Estimation in functional regression for general exponential families}.
\bjournal{The Annals of Statistics}
\bpages{2421--2451}.
\end{barticle}
\endbibitem

\bibitem[\protect\citeauthoryear{Ferraty and
  Vieu}{2006}]{ferraty2006nonparametric}
\begin{bbook}[author]
\bauthor{\bsnm{Ferraty},~\bfnm{Fr{\'e}d{\'e}ric}\binits{F.}} \AND
  \bauthor{\bsnm{Vieu},~\bfnm{Philippe}\binits{P.}}
(\byear{2006}).
\btitle{Nonparametric Functional Data Analysis: Theory and Practice},
\bedition{1} ed.
\bpublisher{Springer-Verlag}, \baddress{New York}.
\end{bbook}
\endbibitem

\bibitem[\protect\citeauthoryear{Giné and Nickl}{2015}]{ginebook}
\begin{bbook}[author]
\bauthor{\bsnm{Giné},~\bfnm{Evarist}\binits{E.}} \AND
  \bauthor{\bsnm{Nickl},~\bfnm{Richard}\binits{R.}}
(\byear{2015}).
\btitle{Mathematical Foundations of Infinite-Dimensional Statistical Models},
\bedition{1} ed.
\bseries{Cambridge Series in Statistical and Probabilistic Mathematics}.
\bpublisher{Cambridge University Press}, \baddress{New York}.
\end{bbook}
\endbibitem

\bibitem[\protect\citeauthoryear{Hall and Horowitz}{2007}]{hall2007methodology}
\begin{barticle}[author]
\bauthor{\bsnm{Hall},~\bfnm{Peter}\binits{P.}} \AND
  \bauthor{\bsnm{Horowitz},~\bfnm{Joel~L.}\binits{J.~L.}}
(\byear{2007}).
\btitle{{Methodology and convergence rates for functional linear regression}}.
\bjournal{The Annals of Statistics}
\bvolume{35}
\bpages{70 -- 91}.
\end{barticle}
\endbibitem

\bibitem[\protect\citeauthoryear{Hall, M{\"u}ller and
  Wang}{2006}]{hall2006properties}
\begin{barticle}[author]
\bauthor{\bsnm{Hall},~\bfnm{Peter}\binits{P.}},
  \bauthor{\bsnm{M{\"u}ller},~\bfnm{Hans-Georg}\binits{H.-G.}} \AND
  \bauthor{\bsnm{Wang},~\bfnm{Jane-Ling}\binits{J.-L.}}
(\byear{2006}).
\btitle{{Properties of principal component methods for functional and
  longitudinal data analysis}}.
\bjournal{The Annals of Statistics}
\bvolume{34}
\bpages{1493 -- 1517}.
\end{barticle}
\endbibitem

\bibitem[\protect\citeauthoryear{He, Huang and He}{2024}]{he2024penalized}
\begin{barticle}[author]
\bauthor{\bsnm{He},~\bfnm{Shiyuan}\binits{S.}},
  \bauthor{\bsnm{Huang},~\bfnm{Jianhua~Z}\binits{J.~Z.}} \AND
  \bauthor{\bsnm{He},~\bfnm{Kejun}\binits{K.}}
(\byear{2024}).
\btitle{Penalized spline estimation of principal components for sparse
  functional data: Rates of convergence}.
\bjournal{Bernoulli}
\bvolume{30}
\bpages{2795--2820}.
\end{barticle}
\endbibitem

\bibitem[\protect\citeauthoryear{Hsing and Eubank}{2015}]{hsing2015theoretical}
\begin{bbook}[author]
\bauthor{\bsnm{Hsing},~\bfnm{Tailen}\binits{T.}} \AND
  \bauthor{\bsnm{Eubank},~\bfnm{Randall}\binits{R.}}
(\byear{2015}).
\btitle{Theoretical foundations of functional data analysis, with an
  introduction to linear operators},
\bedition{1} ed.
\bpublisher{John Wiley \& Sons}, \baddress{Chichester, West Sussex}.
\end{bbook}
\endbibitem

\bibitem[\protect\citeauthoryear{Huang}{2003}]{huang2003local}
\begin{barticle}[author]
\bauthor{\bsnm{Huang},~\bfnm{Jianhua~Z}\binits{J.~Z.}}
(\byear{2003}).
\btitle{Local asymptotics for polynomial spline regression}.
\bjournal{The Annals of Statistics}
\bvolume{31}
\bpages{1600--1635}.
\end{barticle}
\endbibitem

\bibitem[\protect\citeauthoryear{James, Hastie and
  Sugar}{2000}]{james2000principal}
\begin{barticle}[author]
\bauthor{\bsnm{James},~\bfnm{Gareth~M}\binits{G.~M.}},
  \bauthor{\bsnm{Hastie},~\bfnm{Trevor~J}\binits{T.~J.}} \AND
  \bauthor{\bsnm{Sugar},~\bfnm{Catherine~A}\binits{C.~A.}}
(\byear{2000}).
\btitle{Principal component models for sparse functional data}.
\bjournal{Biometrika}
\bvolume{87}
\bpages{587--602}.
\end{barticle}
\endbibitem

\bibitem[\protect\citeauthoryear{James, Wang and
  Zhu}{2009}]{james2009functional}
\begin{barticle}[author]
\bauthor{\bsnm{James},~\bfnm{Gareth~M.}\binits{G.~M.}},
  \bauthor{\bsnm{Wang},~\bfnm{Jing}\binits{J.}} \AND
  \bauthor{\bsnm{Zhu},~\bfnm{Ji}\binits{J.}}
(\byear{2009}).
\btitle{{Functional linear regression that’s interpretable}}.
\bjournal{The Annals of Statistics}
\bvolume{37}
\bpages{2083 -- 2108}.
\end{barticle}
\endbibitem

\bibitem[\protect\citeauthoryear{Jirak and Wahl}{2020}]{JIRAK2020}
\begin{barticle}[author]
\bauthor{\bsnm{Jirak},~\bfnm{Moritz}\binits{M.}} \AND
  \bauthor{\bsnm{Wahl},~\bfnm{Martin}\binits{M.}}
(\byear{2020}).
\btitle{Perturbation bounds for eigenspaces under a relative gap condition}.
\bjournal{Proceedings of the American Mathematical Society}
\bvolume{148}
\bpages{479 -- 494}.
\end{barticle}
\endbibitem

\bibitem[\protect\citeauthoryear{Jirak and Wahl}{2023}]{JIRAK2023108808}
\begin{barticle}[author]
\bauthor{\bsnm{Jirak},~\bfnm{Moritz}\binits{M.}} \AND
  \bauthor{\bsnm{Wahl},~\bfnm{Martin}\binits{M.}}
(\byear{2023}).
\btitle{Relative perturbation bounds with applications to empirical covariance
  operators}.
\bjournal{Advances in Mathematics}
\bvolume{412}
\bpages{108808}.
\end{barticle}
\endbibitem

\bibitem[\protect\citeauthoryear{Kokoszka and
  Reimherr}{2017}]{kokoszka2017introduction}
\begin{bbook}[author]
\bauthor{\bsnm{Kokoszka},~\bfnm{Piotr}\binits{P.}} \AND
  \bauthor{\bsnm{Reimherr},~\bfnm{Matthew}\binits{M.}}
(\byear{2017}).
\btitle{Introduction to Functional Data Analysis},
\bedition{1} ed.
\bpublisher{Chapman and Hall/CRC}, \baddress{Boca Raton}.
\end{bbook}
\endbibitem

\bibitem[\protect\citeauthoryear{Li and Hsing}{2010}]{li2010uniform}
\begin{barticle}[author]
\bauthor{\bsnm{Li},~\bfnm{Yehua}\binits{Y.}} \AND
  \bauthor{\bsnm{Hsing},~\bfnm{Tailen}\binits{T.}}
(\byear{2010}).
\btitle{Uniform convergence rates for nonparametric regression and principal
  component analysis in functional/longitudinal data}.
\bjournal{The Annals of Statistics}
\bvolume{38}
\bpages{3321--3351}.
\end{barticle}
\endbibitem

\bibitem[\protect\citeauthoryear{Paul and Peng}{2009}]{paul2009consistency}
\begin{barticle}[author]
\bauthor{\bsnm{Paul},~\bfnm{Debashis}\binits{D.}} \AND
  \bauthor{\bsnm{Peng},~\bfnm{Jie}\binits{J.}}
(\byear{2009}).
\btitle{{Consistency of restricted maximum likelihood estimators of principal
  components}}.
\bjournal{The Annals of Statistics}
\bvolume{37}
\bpages{1229 -- 1271}.
\end{barticle}
\endbibitem

\bibitem[\protect\citeauthoryear{Peng and Paul}{2009}]{peng2009geometric}
\begin{barticle}[author]
\bauthor{\bsnm{Peng},~\bfnm{Jie}\binits{J.}} \AND
  \bauthor{\bsnm{Paul},~\bfnm{Debashis}\binits{D.}}
(\byear{2009}).
\btitle{A geometric approach to maximum likelihood estimation of the functional
  principal components from sparse longitudinal data}.
\bjournal{Journal of Computational and Graphical Statistics}
\bvolume{18}
\bpages{995--1015}.
\end{barticle}
\endbibitem

\bibitem[\protect\citeauthoryear{Ramsay and Silverman}{2005}]{ramsay2005}
\begin{bbook}[author]
\bauthor{\bsnm{Ramsay},~\bfnm{Jim}\binits{J.}} \AND
  \bauthor{\bsnm{Silverman},~\bfnm{Bernard}\binits{B.}}
(\byear{2005}).
\btitle{Functional Data Analysis},
\bedition{1} ed.
\bseries{Springer Series in Statistics}.
\bpublisher{Springer New York}.
\end{bbook}
\endbibitem

\bibitem[\protect\citeauthoryear{Stone}{1994}]{stone1994use}
\begin{barticle}[author]
\bauthor{\bsnm{Stone},~\bfnm{Charles~J}\binits{C.~J.}}
(\byear{1994}).
\btitle{The use of polynomial splines and their tensor products in multivariate
  function estimation}.
\bjournal{The Annals of Statistics}
\bpages{118--171}.
\end{barticle}
\endbibitem

\bibitem[\protect\citeauthoryear{Wahl}{2022}]{wahl2022lower}
\begin{binproceedings}[author]
\bauthor{\bsnm{Wahl},~\bfnm{Martin}\binits{M.}}
(\byear{2022}).
\btitle{Lower bounds for invariant statistical models with applications to
  principal component analysis}.
In \bbooktitle{Annales de l'Institut Henri Poincare (B) Probabilites et
  statistiques}
\bvolume{58}
\bpages{1565--1589}.
\bpublisher{Institut Henri Poincar{\'e}}.
\end{binproceedings}
\endbibitem

\bibitem[\protect\citeauthoryear{Wang, Chiou and
  M{\"u}ller}{2016}]{wang2016functional}
\begin{barticle}[author]
\bauthor{\bsnm{Wang},~\bfnm{Jane-Ling}\binits{J.-L.}},
  \bauthor{\bsnm{Chiou},~\bfnm{Jeng-Min}\binits{J.-M.}} \AND
  \bauthor{\bsnm{M{\"u}ller},~\bfnm{Hans-Georg}\binits{H.-G.}}
(\byear{2016}).
\btitle{Functional Data Analysis}.
\bjournal{Annual Review of Statistics and Its Application}
\bvolume{3}
\bpages{257--295}.
\end{barticle}
\endbibitem

\bibitem[\protect\citeauthoryear{Wang, Wong and Zhang}{2022}]{wang2022low}
\begin{barticle}[author]
\bauthor{\bsnm{Wang},~\bfnm{Jiayi}\binits{J.}},
  \bauthor{\bsnm{Wong},~\bfnm{Raymond~KW}\binits{R.~K.}} \AND
  \bauthor{\bsnm{Zhang},~\bfnm{Xiaoke}\binits{X.}}
(\byear{2022}).
\btitle{Low-rank covariance function estimation for multidimensional functional
  data}.
\bjournal{Journal of the American Statistical Association}
\bvolume{117}
\bpages{809--822}.
\end{barticle}
\endbibitem

\bibitem[\protect\citeauthoryear{Xiao}{2020}]{xiao2020asymptotic}
\begin{barticle}[author]
\bauthor{\bsnm{Xiao},~\bfnm{Luo}\binits{L.}}
(\byear{2020}).
\btitle{Asymptotic properties of penalized splines for functional data}.
\bjournal{Bernoulli}
\bvolume{26}
\bpages{2847--2875}.
\end{barticle}
\endbibitem

\bibitem[\protect\citeauthoryear{Xiao et~al.}{2018}]{xiao2018fast}
\begin{barticle}[author]
\bauthor{\bsnm{Xiao},~\bfnm{Luo}\binits{L.}},
  \bauthor{\bsnm{Li},~\bfnm{Cai}\binits{C.}},
  \bauthor{\bsnm{Checkley},~\bfnm{William}\binits{W.}} \AND
  \bauthor{\bsnm{Crainiceanu},~\bfnm{Ciprian}\binits{C.}}
(\byear{2018}).
\btitle{Fast covariance estimation for sparse functional data}.
\bjournal{Statistics and computing}
\bvolume{28}
\bpages{511--522}.
\end{barticle}
\endbibitem

\bibitem[\protect\citeauthoryear{Yan, Yao and Zhou}{2025}]{yan2025drrm}
\begin{barticle}[author]
\bauthor{\bsnm{Yan},~\bfnm{Shunxing}\binits{S.}},
  \bauthor{\bsnm{Yao},~\bfnm{Fang}\binits{F.}} \AND
  \bauthor{\bsnm{Zhou},~\bfnm{Hang}\binits{H.}}
(\byear{2025}).
\btitle{Deep regression for repeated measurements}.
\bjournal{Journal of the American Statistical Association}
\bvolume{120}
\bpages{2461--2472}.
\end{barticle}
\endbibitem

\bibitem[\protect\citeauthoryear{Yan and Yao}{2026a}]{yanflr}
\begin{barticle}[author]
\bauthor{\bsnm{Yan},~\bfnm{Shunxing}\binits{S.}} \AND
  \bauthor{\bsnm{Yao},~\bfnm{Fang}\binits{F.}}
(\byear{2026}a).
\btitle{Functional linear regression from sparse to dense designs: a
  pooling-ridge method and minimax optimality}.
\bjournal{The Annals of Statistics}.
\bnote{Forthcoming. \url{https://arxiv.org/abs/2608.25468}}.
\end{barticle}
\endbibitem

\bibitem[\protect\citeauthoryear{Yan and Yao}{2026b}]{supplyment}
\begin{barticle}[author]
\bauthor{\bsnm{Yan},~\bfnm{Shunxing}\binits{S.}} \AND
  \bauthor{\bsnm{Yao},~\bfnm{Fang}\binits{F.}}
(\byear{2026}b).
\btitle{Supplement to ``Functional Principal Component Analysis: Minimax
  Optimality from Sparse to Dense Designs and Beyond''}.
\end{barticle}
\endbibitem

\bibitem[\protect\citeauthoryear{Yao, M{\"u}ller and Wang}{2005a}]{yao2005a}
\begin{barticle}[author]
\bauthor{\bsnm{Yao},~\bfnm{Fang}\binits{F.}},
  \bauthor{\bsnm{M{\"u}ller},~\bfnm{Hans-Georg}\binits{H.-G.}} \AND
  \bauthor{\bsnm{Wang},~\bfnm{Jane-Ling}\binits{J.-L.}}
(\byear{2005}a).
\btitle{Functional data analysis for sparse longitudinal data}.
\bjournal{Journal of the American statistical association}
\bvolume{100}
\bpages{577--590}.
\end{barticle}
\endbibitem

\bibitem[\protect\citeauthoryear{Yao, M{\"u}ller and Wang}{2005b}]{yao2005b}
\begin{barticle}[author]
\bauthor{\bsnm{Yao},~\bfnm{Fang}\binits{F.}},
  \bauthor{\bsnm{M{\"u}ller},~\bfnm{Hans-Georg}\binits{H.-G.}} \AND
  \bauthor{\bsnm{Wang},~\bfnm{Jane-Ling}\binits{J.-L.}}
(\byear{2005}b).
\btitle{{Functional linear regression analysis for longitudinal data}}.
\bjournal{The Annals of Statistics}
\bvolume{33}
\bpages{2873 -- 2903}.
\end{barticle}
\endbibitem

\bibitem[\protect\citeauthoryear{Yuan and Cai}{2010}]{yuan2010reproducing}
\begin{barticle}[author]
\bauthor{\bsnm{Yuan},~\bfnm{Ming}\binits{M.}} \AND
  \bauthor{\bsnm{Cai},~\bfnm{T.~Tony}\binits{T.~T.}}
(\byear{2010}).
\btitle{{A reproducing kernel Hilbert space approach to functional linear
  regression}}.
\bjournal{The Annals of Statistics}
\bvolume{38}
\bpages{3412 -- 3444}.
\end{barticle}
\endbibitem

\bibitem[\protect\citeauthoryear{Zhang and Wang}{2016}]{zhang2016sparse}
\begin{barticle}[author]
\bauthor{\bsnm{Zhang},~\bfnm{Xiaoke}\binits{X.}} \AND
  \bauthor{\bsnm{Wang},~\bfnm{Jane-Ling}\binits{J.-L.}}
(\byear{2016}).
\btitle{From sparse to dense functional data and beyond}.
\bjournal{The Annals of Statistics}
\bvolume{44}
\bpages{2281--2321}.
\end{barticle}
\endbibitem

\bibitem[\protect\citeauthoryear{Zhou, Wei and
  Yao}{2025}]{zhou2024theoryfunctionalprincipalcomponent}
\begin{barticle}[author]
\bauthor{\bsnm{Zhou},~\bfnm{Hang}\binits{H.}},
  \bauthor{\bsnm{Wei},~\bfnm{Dongyi}\binits{D.}} \AND
  \bauthor{\bsnm{Yao},~\bfnm{Fang}\binits{F.}}
(\byear{2025}).
\btitle{{Theory of functional principal component analysis for discretely
  observed data}}.
\bjournal{The Annals of Statistics}
\bvolume{53}
\bpages{2103 -- 2127}.
\end{barticle}
\endbibitem

\bibitem[\protect\citeauthoryear{Zhou, Yao and
  Zhang}{2023}]{zhou2023functional}
\begin{barticle}[author]
\bauthor{\bsnm{Zhou},~\bfnm{Hang}\binits{H.}},
  \bauthor{\bsnm{Yao},~\bfnm{Fang}\binits{F.}} \AND
  \bauthor{\bsnm{Zhang},~\bfnm{Huiming}\binits{H.}}
(\byear{2023}).
\btitle{Functional linear regression for discretely observed data: from ideal
  to reality}.
\bjournal{Biometrika}
\bvolume{110}
\bpages{381--393}.
\end{barticle}
\endbibitem

\end{thebibliography}

\end{document}